\documentclass[11pt]{article}

\usepackage{amsmath,amsfonts}

\newcommand{\be}{\begin}
\newcommand{\Ref}[1]{(\ref{#1})}

\newcommand{\cm}{\mathcal M}

\newcommand{\cb}{\mathcal B}
\newcommand{\cc}{\mathcal C}

\newcommand{\cl}{\mathcal L}
\newcommand\ct{\mathcal T}

\def\fl{\mathcal R}

\newcommand{\al}{\alpha}
\newcommand{\del}{\delta}
\newcommand{\eps}{\epsilon}
\def\ga{\gamma}
\def\Ga{\Gamma}
\def\ka{\kappa}
\def\lla{\lambda}
\def\La{\Lambda}
\def\om{\omega}

\newcommand{\si}{\sigma}
\def\Si{\Sigma}

\newcommand{\co}{\mathbb C}

\newcommand{\e}{{\mathbb E}}
\def\f{{\mathbb F}}

\def\bk{\mathbb K}
\def\bl{\mathbb L}
\def\bs{\mathbb S}
\newcommand\bt{{\mathbb T}}

\def\q{{\mathbb Q}}
\newcommand{\R}{\mathbb R}

\newcommand{\z}{\mathbb Z}

\def\cp#1{{\mathbb{CP}^{#1}}}

\newcommand{\uset}{\underset}
\newcommand{\oset}{\overset}
\newcommand{\uline}{\underline}
\newcommand{\oline}{\overline}

\newcommand{\la}{\langle}
\newcommand{\ra}{\rangle}
\newcommand{\st}{\,:\,}
\newcommand{\ti}{\tilde}
\newcommand{\setmin}{\!\setminus\!}
\newcommand{\prtl}{\partial}
\def\square{\kern20pt{\vbox{\hrule height.4pt
        \hbox{\vrule width.4pt height 6pt\kern6pt
                \vrule width.4pt}
        \hrule height.4pt}}}
\newcommand{\ry}{\R\times Y}
\newcommand{\rpy}{\R_+\times Y}
\newcommand\rmy{\R_-\times Y}

\def\im{\text{\rm im}}

\def\coker{\text{\rm coker}}

\def\inv{^{-1}}

\def\proof{{\em Proof.}\ }

\def\SW{\text{SW}}

\def\hf{\text{HF}}
\def\rhf{\widehat{\text{HF}}\vphantom{F}}
\def\cf{\text{CF}}
\def\ohf{\overline{\text{HF}}\vphantom{\text{F}}}
\def\uhf{\underline{\text{HF}}\vphantom{\text{F}}}

\def\U#1{\text{U}(#1)}

\def\ind{\text{ind}}

\def\spc{\text{spin}^c}

\def\cd{{\vartheta}}
\def\ucd{\uline{\cd}}
\def\uucd{\uline{\uline{\cd}}}
\def\csd#1{\vartheta_{#1}}
\def\mm{\mathfrak m}
\def\prt{\mathfrak p}
\def\Prt{\mathfrak P}
\def\s{\mathfrak s}

\newcommand\llw[3]{{L^{#1,#2}_#3}}
\newcommand\lw[2]{{L^{#1,#2}}}

\def\thet#1{\theta_{m_#1}}

\def\hm{\hat M}
\def\hw{\hat W}
\newcommand\hz{\hat Z}
\newcommand\hatF{\hat F}
\def\rms{\R_-\times S^3}
\def\id{\text{id}}

\def\utr{U_T}
\def\gtr{g_{T,r}}
\def\gcyl{g_u^c}
\def\gy{g}
\def\gtru{g_{T,r,u}}
\def\dtr{D_{T,r}}
\def\dtru{D_{T,r,u}}
\def\tgtru{\ti g_{T,r,u}}

\def\opsi{\oline\psi}
\def\upsi{\uline\psi}

\def\rpsi{\hat\psi}
\def\Ltr{L}
\def\tr{\text{tr}}

\def\ww#1#2{W_{#1,#2}}
\def\zinf{_{0,\infty}}
\def\jinf{_{j,\infty}}
\def\jn{_{j,n}}
\newcommand\tu{^{(T)}}
\newcommand\xt{{X\tu}}
\def\xtn{X^{(T(n))}}

\def\ocan{o_{\text{can}}}
\newcommand\chm{\check M}
\def\chv#1{\check V_#1}
\def\hut{\hat U_\tau}
\def\ared{A_{\text{red}}}
\def\bred{B_{\text{red}}}
\def\omred{\omega_{\text{red}}}
\newcommand\rr{_o}
\def\uka{\uline\ka}
\def\smo{S^3_{-1}}
\def\sze{S^3_0}
\def\ip{\uline p}

\begin{document}

\title{Monopole Floer homology\\for rational homology $3$--spheres}

\author{Kim A.\ Fr\o yshov}


\date{}

\maketitle

\be{abstract}
We give a new construction of monopole Floer homology 
for $\spc$ rational homology $3$--spheres.
As applications we define two invariants of certain $4$--manifolds with
$b_1=1$ and $b^+=0$.
\end{abstract}

\bibliographystyle{plain}

\tableofcontents

\section{Introduction}

In their recent book \cite{KM5} Kronheimer and Mrowka
constructed monopole Floer homology
groups for arbitrary closed oriented $3$--manifolds by applying a new kind
of Morse homology to
certain blown-up configuration spaces. In this paper we will develop the 
theory from a different angle in the case of an oriented
rational homology $3$--sphere $Y$ equipped with 
a $\spc$ structure. Our starting-point will be the
``irreducible'' Floer cohomology $\hf^*(Y,m)$ which is the result of
adapting Floer's original definition \cite{F1} to monopoles. Here 
the parameter $m$, which runs through a set $\mm(Y)$ of the form
$m_0+\z$ with $m_0\in\q$,
indicates which chamber for the metric and perturbation is being used.
(Because of our choice of grading conventions we find it more natural to
work with Floer {\em co}homology.)

By taking suitable limits of $\hf^*(Y,m)$ as
$m\to\pm\infty$ we obtain invariants $\ohf^*(Y)$, $\uhf^*(Y)$ 
of the $\spc$ manifold $Y$. These ``equivariant'' Floer cohomology groups are
modules over a polynomial ring and are related by a long exact sequence
(which we call the {\em fundamental sequence})
involving a third module which is essentially the module of Laurent polynomials.
Exact sequences of this kind are well known from
Ozsv\'ath and Szab\'o's Heegaard Floer homology \cite{OS3} and also appear
in \cite{KM5}. In fact, it seems likely that our
equivariant groups are isomorphic to the ($*$-versions of the)
monopole Floer groups of Kronheimer and Mrowka (and their Heegaard cousins)
when differences in grading
conventions are taken into account.

A perhaps unusual feature in our model is that the fundamental sequence
is not constructed from a short exact sequence of cochain complexes.
Exactness is established by geometric means,
without recourse to homological algebra.

In the case of coefficients in a field $\f$ of characteristic $p$
(where we allow $p=0$)
we give a precise description of how $\hf^*(Y,m;\f)$ depends on $m$, in terms
of an invariant $h_p(Y)\in\mm(Y)$. For a fixed chamber $m$ the integer
$m-h_p(Y)$ measures interaction between irreducible critical points and 
the reducible critical point of the Chern--Simons--Dirac functional,
much in the same way as the instanton $h$--invariant (see \cite{Fr3}).
The invariant $h_p(Y)$ can also be read off from the fundamental sequence.

Our approach to understanding the relationship between the groups
$\hf^*(Y,m)$ for different $m$ involves a trick to avoid obstructed
gluing problems. The idea is to study moduli spaces not only over the
cylinder $\ry$ but also over a manifold obtained by removing one point
from $\ry$, thereby creating an additional end $\R_-\times S^3$. By placing
a suitable metric and perturbation on this new end we are able to
correct the index of the Dirac operator. In this way we construct certain
chain homomorphisms which are also essential to
the proof of the exactness of the fundamental sequence.

In a slightly different direction we give some applications of Floer cohomology
to $\spc$ smooth closed $4$--manifolds $X$
with $b_1(X)=1$ which contain a non-separating smoothly embedded
rational homology $3$--sphere $Y$.
In the case $b^+(X)>1$ we express the Seiberg--Witten
invariant of $X$ as the Lefschetz number of a certain endomorphism
of the {\em reduced} Floer cohomology of $Y$. If $b^+(X)=0$ 
then, for certain $\spc$ structures on $X$, this Lefschetz number yields
an invariant of $X$, i.e.\ is independent of $Y$. If in addition $Y$
is an integral homology sphere or $b_2(X)=0$
then $h_p(Y)$ is an invariant of $(X,e)$,
where $e$ denotes the generator of $H_3(X;\z)$ represented by $Y$.

Part of this paper was written in 2001-2 while the author was 
a visitor at the Institut des Hautes \'Etudes Scientifiques, to which
the author is grateful for its hospitality.
This work was also partially supported by NSF grant DMS-9971731,
a grant from the
DFG (German Research Foundation) as well as by the 
CRC 701 at the University of Bielefeld.

{\bf Acknowledgement:} The author wishes to thank an anonymous referee
for several suggestions that helped improve the exposition of the paper.

\section{Main results}

\subsection{Irreducible Floer groups}

By a $\spc$ manifold we shall mean a smooth, oriented manifold with a 
$\spc$ structure. The $\spc$ structure will usually be suppressed from
notation, except in explicit computations. Let $Y$ be a $\spc$ rational
homology $3$--sphere. 
There is a classical invariant $\mm(Y)\in\q/\z$ which may be characterized
as follows: If $X$ is any compact $\spc$ $4$--manifold with
$\prtl X=Y$ (as $\spc$ manifolds) then 
\[\mm(Y)\equiv\frac18(\si(X)-c_1(\cl_X)^2)\mod\z.\]
Here $\cl_X$ is the determinant line bundle of the $\spc$ structure,
and $\si(X)$ is the signature of $X$. We will think of $\mm(Y)$ as
a set of rational numbers of the form $m_0+\z$, $m_0\in\q$. 
Note that $\mm(Y)=\z$ if $Y$ is an integral homology sphere.

For any $m\in\mm(Y)$ we will define
an {\em irreducible} Floer cohomology group $\hf^*(Y,m)$ which is
a finitely generated, graded Abelian group, where the grading runs through
$2m+\z$. This is the
monopole analogue of the instanton homology groups introduced by Floer
\cite{F1}. We will use the set-up of \cite{Fr13}, see in particular
Section~1.2 and Chapter~3 of that book.

Choose a Riemannian metric $g$ and $1$--form $\nu$ on $Y$, and define the
Chern--Simons--Dirac functional $\cd$ as in \cite[Section~3.2]{Fr13}
with $\eta=d\nu$. Let $\fl=\fl_{(g,\nu)}=\fl_Y$ be the set of 
gauge equivalence classes of critical points of $\cd$. There is a unique
reducible point in $\fl$, which we denote by $\theta$.
Set $\fl^*:=\fl\setmin\{\theta\}$. In Section~\ref{sec:modsp-ind}
we will associate to any non-degenerate
$\al\in\fl$ an index (or degree) in $2\mm(Y)+\z$.
For generic $\nu$, all critical points of $\csd{}$ are non-degenerate,
in which case $\fl$ is a finite set.

We can then form a cochain complex
$\cf^*$ where $\cf^q$ is essentially the free Abelian group generated by the 
elements of $\fl^*$ of index $q$, and the differential 
is defined in terms of Seiberg--Witten moduli spaces over the cylinder $\ry$.
The cohomology group of $\cf^*$ turns out to depend only on the $\spc$
manifold $Y$ and $m:=\frac12\ind(\theta)$.
This cohomology group is what we denote by
$\hf^*(Y,m)$. It is defined for all $m\in\mm(Y)$ because every
such $m$ can be realized for some pair $(g,\nu)$, as we will show in
Section~\ref{sec:chambers}.

Let $G$ be any Abelian group and $\hf^*(Y,m;G)$ the irreducible
Floer cohomology group with coefficients in $G$.
If $m_1\le m_2$ then there is a canonical homomorphism
\be{equation*}
J_{m_1}^{m_2}:\hf^*(Y,m_1;G)\to\hf^*(Y,m_2;G)
\end{equation*}
which is an isomorphism in degree $q$ unless $2m_1\le q<2m_2$.
These maps obey the composition law
$J_{m_2}^{m_3}\circ J_{m_1}^{m_2}=J_{m_1}^{m_3}$. 

\be{thm}\label{hp-char}
Let $\f$ be a field of characteristic $p$. Then
there is an element $h_p(Y)\in\mm(Y)$ such that
for Floer cohomology with coefficients in $\f$ the following hold:
\be{description}
\item[(i)]Let $m_1\le m_2\le h_p(Y)$. Then $J^{m_2}_{m_1}$
is surjective, and its kernel has dimension~$1$ in degree $2m_1+2k-1$ for
$k=1,\dots,m_2-m_1$ and is zero in all other degrees.
\item[(ii)]Let $h_p(Y)\le m_1\le m_2$. Then $J^{m_2}_{m_1}$ is injective,
and its cokernel has dimension~$1$ in degree $2m_2-2k$ for 
$k=1,\dots,m_2-m_1$ and is zero in all other degrees.
\end{description}
\end{thm}

A precise description of the kernel and cokernel of $J_{m_1}^{m_2}$ is 
provided by Proposition~\ref{prop:varchamb} below. Some related material
can be found in \cite[Section~7]{Marcolli-Wang1}.

It follows from Proposition~\ref{prop:Jtorsion} below that if
$\hf^*(Y,h_0(Y);\z)$ is torsion-free then $\hf^*(Y,m;\z)$ is torsion-free 
for all $m$. Now, if $Y$ admits a metric of positive scalar curvature,
then $\hf^*(Y,h_0(Y);\z)=0$.
Applying Theorem~\ref{hp-char} with rational coefficients then yields
a complete description of $\hf^*(Y,m;\z)$ for all $m$.

\subsection{Equivariant Floer groups}
\label{subsec:equiv-floergr}

We define the 
equivariant Floer cohomology groups $\ohf^q(Y;G)$ and $\uhf^q(Y;G)$
to be the direct and inverse limits, respectively, of
the sequence of homomorphisms
\[\dots\to\hf^q(Y,m-1;G)\to\hf^q(Y,m;G)\to\hf^q(Y,m+1;G)\to\dots.\]
(Section~\ref{sec:fd-analogue} will provide some justification for the term
``equivariant''.)
Clearly, there are natural homomorphisms
\[\uhf^q(Y;G)\to\hf^q(Y,m;G)\to\ohf^q(Y;G)\]
for any $q,m$. The first of these maps is an isomorphism when $2m\le q$,
the second one when $q\le2m-1$. Hence 
$\uhf^q(Y;G)$ vanishes for $q\gg0$, and 
$\ohf^q(Y;G)$ vanishes for $q\ll0$.
We define the homology groups $\ohf_q(Y;G)$ and $\uhf_q(Y;G)$ in a similar
fashion as limits of the groups $\hf_q(Y,m;G)$ as $m\to\mp\infty$.

The equivariant groups $\uhf^*(Y;G)$ and $\ohf^*(Y;G)$ both
come equipped with a natural degree~$2$ endomorphism
$u$. There is also a natural homomorphism of $\z[u]$--modules
\[J=J_Y:\uhf^*(Y;G)\to\ohf^*(Y;G)\]
whose image we denote by $\rhf^*(Y;G)$ and call the {\em reduced} Floer 
cohomology of $Y$.

The two kinds of equivariant Floer groups are related by a long exact 
sequence (the {\em fundamental sequence})
which we now describe. This sequence involves the graded
$\z[u]$--module $P^*(Y)$ of ``generalized polynomials''
\[\sum_{m\in\mm(Y)}a_mx^m\]
where the coefficients $a_m$ are integers and only finitely many of them
are non-zero. The grading and module structure are given by
\[\deg(x^m)=2m,\qquad u\cdot x^m=x^{m+1}.\]
If $Y$ is an integral homology sphere then
$\mm(Y)=\z$, so in this case $P^*(Y)$ is the module of Laurent polynomials.
There is a canonical long exact sequence of $\z[u]$--homomorphisms
\be{equation}\label{eqn:exact-seq}
\cdots\oset{D'}\longrightarrow\uhf^*(Y;G)\oset J\longrightarrow\ohf^*(Y;G)
\oset D\longrightarrow P^*(Y)\otimes G\oset{D'}\longrightarrow\uhf^{*+1}(Y;G)
\oset J\longrightarrow\cdots
\end{equation}
where $J,D$ have degree $0$ and $D'$ has degree $1$. The maps
$D,D'$, which are defined at the end of Section~\ref{sec:diff-chambers},
in some sense measure interaction between reducible and irreducible
critical points of the Chern--Simons--Dirac functional. In our set-up this
sequence is not constructed from a short exact sequence of cochain complexes.
Instead, exactness is proved by means of an ``almost inverse'' of the map
$J_m^{m+1}$.

\be{thm}\label{thm:h-D}
If $\f$ is any field of characteristic $p$ then
$2h_p(Y)$ is the lowest degree in which $D:\ohf^*(Y;\f)\to P^*(Y)\otimes\f$
is non-zero.
\end{thm}

Thus, with field coefficients the irreducible Floer groups contain no
information
that cannot be read off from the sequence \Ref{eqn:exact-seq}. (With
integral coefficients the situation might be different.)

Recall that reversing the orientation of $Y$ reverses the sign of the
Chern--Simons--Dirac functional. This gives rise to canonical
``Poincar\'e duality'' isomorphisms
\be{gather*}
\hf^*(-Y,-m;G)=\hf_{-1-*}(Y,m;G),\\
\ohf^*(-Y;G)=\uhf_{-1-*}(Y;G),\quad
\uhf^*(-Y;G)=\ohf_{-1-*}(Y;G),
\end{gather*}
the latter two being isomorphisms of $\z[u]$--modules. 

\subsection{The $h$--invariant}

In this section $p$ will be fixed and we mostly write
$h$ instead of $h_p$.

\be{thm}\label{thm:h-add}For all $\spc$ rational homology $3$--spheres
$Y_1,Y_2$ one has
\[h(Y_1\#Y_2)=h(Y_1)+h(Y_2).\]
\end{thm}

\be{thm}\label{main-ineq}
Let $W$ be a smooth, compact $\spc$ $4$--manifold whose
boundary is a disjoint union of
rational homology spheres $Y_1,\dots,Y_r$. If $W$ has negative definite
intersection form then
\[-\sum_{j=1}^rh(Y_j)\ge\frac18(b_2(W)+c_1(\cl_W)^2).\]
\end{thm}

Note that if $d$ denotes the correction term in Heegaard Floer homology
then the above two theorems hold with $-d/2$ in place of $h$,
see \cite{OS2}.

If each $Y_j$ is an integral homology sphere then by a theorem
of Elkies~\cite{Elkies1} there is a $\spc$ structure on $W$ for which
the right hand side of the inequality in Theorem~\ref{main-ineq} is
non-negative,
and positive if the intersection form is not diagonalizable over the integers.

In the case of binary polyhedral spaces, \cite[Proposition~8]{Fr1}
holds with $-8h$ in place of $\ga$. In particular,
for the Poincar\'e homology sphere oriented as the boundary
of the negative definite $E_8$--manifold one has
\[h(\Si(2,3,5))=-1.\]

For simple lens spaces $L(q,1)$, $q\ge2$, Theorem~\ref{main-ineq}
and the property $h(-Y)=-h(Y)$ suffice for the computation of $h$.
To see this, let $\Si$ be an internal connected sum of $q$ disjoint
$(-1)$--spheres in the $q$--fold connected sum $W:=-q\cp2$. Then the link
of $\Si$, a copy of $L(q,1)$, separates $W$ into two pieces $W_1,W_2$.
Now observe that every $\spc$ structure on $L(q,1)$ is the restriction
of a $\spc$ structure $\s$ on $W$ with $c_1(\s)^2=-q$.
Applying Theorem~\ref{main-ineq} to each of $W_1,W_2$ then gives
\[h(L(q,1),\s_j)=\frac18\left[\frac{(q-2j)^2}q-1\right]\]
for a suitable labelling $\s_0,\dots,\s_{q-1}$ of the $\spc$ structures 
on $L(q,1)$.

This example shows that if $W$ is a negative definite $\spc$ cobordism 
from one rational homology $3$--sphere $Y_1$ to another one $Y_2$
(i.e.\ $\prtl W=(-Y_1)\cup Y_2$)
then one need not have $h(Y_1)\ge h(Y_2)$ unless $Y_1,Y_2$ are integral
homology spheres.

There is a relationship between the $h$--invariant and Casson's invariant
which we will now describe.
It follows from Theorem~\ref{hp-char} that
$\ti\lla(Y):=\chi(\hf^*(Y,m))-m$ is independent of $m$ and therefore
an invariant of the $\spc$ manifold $Y$.
(Here $\chi$ denotes the Euler characteristic with respect to the
mod~$2$ grading.)  This result was proved with different methods by
Lim \cite{Lim1}. He also showed that, when $Y$ is an integral homology
sphere, $-\ti\lla$ agrees with Casson's invariant $\lla$
(normalized so that $\lla(\Si(2,3,5))=-1$), thereby confirming a 
conjecture of Kronheimer. Since $\hf^*(Y,h_p(Y);\f)$ maps isomorphically
onto $\rhf^*(Y;\f)$, we obtain:

\be{thm}\label{thm:h-casson}
Let $\f$ be a field of characteristic $p$. Then
for every oriented integral homology $3$--sphere $Y$ one has
\[h_p(Y)-\chi(\rhf^*(Y;\f))=\lla(Y),\]
where $\lla$ denotes Casson's invariant.
\end{thm}

Similar results hold for the instanton and Heegaard Floer homologies, see
\cite[Section~8]{Fr3} and \cite[Theorem~1.3]{OS2}.

\subsection{Invariants of $4$--manifolds}
\label{subsec:inv4mfld}

The results in this subsection will only be stated in their simplest form,
ignoring the $1$--dimensional $\mu$--classes, and using integral coefficients
for the Floer groups unless otherwise specified.
Throughout the subsection $Y,Y'$ will denote rational homology $3$--spheres.
All $4$--manifolds will be smooth.

Let $Z$ be a compact, homology oriented $\spc$ $4$--manifold with $\prtl Z=Y$.
Then $Z$ has a ``relative'' Seiberg--Witten invariant
\[\opsi(Z)\in\ohf^{-d}(Y),\]
where
\[d=\frac14\left(c_1(\cl_Z)^2-\si(Z)\right)+b_1(Z)-b^+(Z).\]
If moreover $b^+(Z)>1$ then there is also an invariant
\[\upsi(Z)\in\uhf^{-d}(Y)\]
with $J\upsi(Z)=\opsi(Z)$.

\be{thm}\label{thm:basic}
Let $Z$ be a closed, connected $\spc$ $4$--manifold
which is separated by an embedded rational homology $3$--sphere $Y$:
\[Z=Z_0\cup_YZ_1.\]
Let $Z_0$ and $Z_1$ be homology oriented and let $Z$ be equipped with the
corresponding glued homology orientation as defined in 
\cite[Section~12.6]{Fr13}.
If $b^+(Z_1)>1$ and $\dim M(Z)=2n\ge0$,
then the Seiberg--Witten invariant of $Z$ is
\[\SW(Z)=u^n\opsi(Z_0)\cdot\upsi(Z_1).\]
\end{thm}

We now turn to cobordisms. Let $W$ be a compact, homology oriented
$\spc$ $4$--manifold with $\prtl W=(-Y)\cup Y'$. Then $W$ induces
graded homomorphisms of $\z[u]$--modules
\[\upsi(W):\uhf^*(Y)\to\uhf^{*-d}(Y'),\qquad
\opsi(W):\ohf^*(Y)\to\ohf^{*-d}(Y')\]
which intertwine with the $J$--maps and therefore induce a homomorphism
\[\rpsi(W):\rhf^*(Y)\to\rhf^{*-d}(Y').\]
Here $d$ is defined as above with $W$ in place of $Z$.
In general, the $\opsi$--invariant of
a composite cobordism or of a composite manifold $Z\cup_YW$ is the 
composition (or product) of the $\opsi$--invariants of the two pieces,
provided the homology orientations are related as in Theorem~\ref{thm:basic}.
The same holds for the $\upsi$--invariants as long as these are defined.
If $b^+(W)>0$ then $\upsi(W)$ factors through $J_Y$.

If $b_1(W)=0$ or $b^+(W)>0$ we define the $\z[u]$--homomorphism
\[P(W):P^*(Y_1)\to P^{*-d}(Y_2)\]
as follows: If $b^+(W)>0$ set $P(W):=0$. If $b^+(W)=0=b_1(W)$ set
$P(W)\cdot x^m:=x^{m-d/2}$. Then we have a commutative diagram
\be{equation}\label{diagr:W}
\be{array}{ccccccc}
\uhf^*(Y_1) & \longrightarrow & \ohf^*(Y_1) &
\longrightarrow & P^*(Y_1) & \longrightarrow & \uhf^{*+1}(Y_1) \\
\downarrow & & \downarrow & & \downarrow & & \downarrow \\
\uhf^{*-d}(Y_2) & \longrightarrow & \ohf^{*-d}(Y_2) &
\longrightarrow & P^{*-d}(Y_2) & \longrightarrow
& \uhf^{*+1-d}(Y_2)
\end{array}
\end{equation}
where the vertical maps are the ones induced by $W$ and the horizontal
ones are as in \Ref{eqn:exact-seq}. We expect that there is a similar
commutative diagram when $b^+(W)=0$, $b_1(W)>0$
but have not yet computed the map $P(W)$ in that case. (The left-most square
of the diagram commutes for any $W$, as already mentioned.)

The homomorphisms induced by cobordisms can be defined more generally
for Floer groups with coefficients in any Abelian group, and have similar
properties.

We will now describe an analogue of Theorem~\ref{thm:basic} for a 
non-separating hypersurface $Y$.

If $V=V_0\oplus V_1$ is any finitely generated mod~2 graded Abelian group
then the {\em Lefschetz number} of a degree preserving
endomorphism $f=f_0\oplus f_1$ of $V$ is defined as
\[\Ltr(f):=\tr(f_0)-\tr(f_1).\]

If $Y$ is an oriented compact hypersurface in an oriented manifold $X$ then
we denote by $X_{\|Y}$ the oriented manifold with boundary $(-Y)\cup Y$
obtained by cutting $X$ open along $Y$.

\be{thm}\label{thm:swltr}
Let $X$ be a closed, connected, homology oriented $\spc$ $4$--manifold
with $b^+(X)>1$. Let $Y$ be a non-separating, smoothly embedded oriented
rational homology $3$--sphere in $X$ and set $W=X_{\|Y}$.
If $\dim M(X)=2n\ge0$ then
\be{equation}\label{eqn:sw-ltr}
\SW(X)=\Ltr(u^n\rpsi(W)),
\end{equation}
where $\rpsi(W):\rhf^*(Y)\to\rhf^{*-2n}(Y)$.
\end{thm}

Here, as well as in Theorem~\ref{thm:Lhpwelldef} below,
$W$ should have the homology orientation determined by the
homology orientation of $X$ and the orientation of $Y$ as specified in
\cite[Section~12.5]{Fr13}.

The fact that the right hand side of \Ref{eqn:sw-ltr} is defined whenever 
$\dim M(X)$ is non-negative and even can be exploited to extend the
Seiberg--Witten invariant to a class of $4$--manifolds for which the usual
definition is not available.

\be{thm}\label{thm:Lhpwelldef}
Let $X$ be a closed, connected, homology oriented $\spc$ $4$--manifold
with $b_1(X)=1$, $b^+(X)=0$, and $c_1(\cl_X)^2=\si(X)$, 
so that $\dim M(X)=0$. For $j=0,1$ let $Y_j$ be a non-separating,
smoothly embedded oriented rational homology $3$--sphere in
$X$. Suppose $Y_0$ and $Y_1$ represent the same generator of $H_3(X;\z)$,
and set $W_j=X_{\|Y_j}$. Then
\[\Ltr(\rpsi(W_0))=\Ltr(\rpsi(W_1)).\]
If $Y_0$ and $Y_1$ are in fact integral homology spheres, or if $b_2(X)=0$, then
\[h_p(Y_0)=h_p(Y_1)\]
for all $p$.
\end{thm}

Assuming the appropriate hypersurface $Y\subset X$ exists
(which is not always the
case) we can therefore define invariants
\[L(X):=\Ltr(\rpsi(W)),\qquad h_p(X,e):=h_p(Y),\]
where $W=X_{\|Y}$, and $e$ is the generator of $H_3(X;\z)$ represented by $Y$.
For instance, if $Y$ is any oriented integral homology $3$--sphere and 
$S^1\times Y$ has the homology orientation opposite to the one given by
the canonical generators of $H^0(S^1\times Y)$ and $H^1(S^1\times Y)$ then
by Theorem~\ref{thm:h-casson} one has
\[L(S^1\times Y)=\chi(\rhf^*(Y;\z))=\chi(\rhf^*(Y;\q))=h_0(Y)-\lla(Y).\]
In particular,
\[L(S^1\times S^3)=0,\qquad L(S^1\times\Si(2,3,7))=1.\]
Moreover, the global spherical shell conjecture for minimal class VII
surfaces with $b_2>0$ (known for $b_2=1$, see \cite{Teleman2})
would imply that $L$ is defined and zero for such surfaces with the
canonical $\spc$ structure.

There is an alternative description of $L(X)$
as the number of irreducible monopoles over $X$, counted with sign,
for a suitable
choice of metric and perturbation (using a long neck $[-T,T]\times Y$).
However, it was shown by Okonek--Teleman \cite{okonek-teleman1}
that this number
in general depends on the
chamber of the metric and perturbation.

The reader may wish to compare the definition of $L(X)$ with 
Ruberman--Saveliev's description \cite{rub-sav1}
of the Furuta--Ohta invariant \cite{Fur-Ohta1} in the case of mapping tori.



\subsection{Outline}

Here is an outline of the remainder of this paper.

Section~3 provides a proof of a technical fact which underpins our whole
approach to Floer homology, namely the non-emptiness of every chamber for
the metric and perturbation on a $\spc$ rational homology $3$--sphere.
The proof relies on B\"ar's gluing theorem for spectra of operators of
Dirac type. The proof may safely be skipped on first reading, since 
the techniques involved are rather different from those in the rest of the
paper. 
Section~4 introduces moduli spaces of monopoles and defines the index
of a non-degenerate critical point of the Chern--Simons--Dirac functional.
Section~5 shows that certain properties of reducible monopoles
over $4$--manifolds are stable
under small perturbations.
In Section~6 we specify our orientation conventions
and describe the behaviour of orientations under (un)gluing maps.
In Section~7 the irreducible Floer cohomology $\hf^*(Y,m)$ is defined
and a first comparison of these groups for different $m$ is made.
Section~8 introduces the $u$--map and establishes its compatibility with
cobordism-induced homomorphisms.
The subject of Section~9 is 
the extra structure on the irreducible Floer groups that arises from the
presence of the reducible critical point, and the way this structure
relates to maps induced by cobordisms.
In Section~10 the trick mentioned in the
introduction is used to determine the kernel and cokernel of 
$J_{m_1}^{m_2}$. At the end of the section the maps $D,D'$ are defined
and the exactness of the fundamental sequence ist established.
In Section~11
it is shown how the reduced Floer cohomology may be computed in an
arbitrary chamber for the metric and perturbation. This is then used
to prove the Lefschetz trace formula of Theorem~\ref{thm:swltr}.
Section~12 contains all material on the $h$--invariant. Most of
this section is logically independent from Sections 9--11.
Section~13 is devoted to 
negative definite $4$--manifolds with $b_1=1$.
The final Section~14 presents a finite-dimensional analogue of the construction
of the equivariant group $\ohf^*(Y)$ in order to justify the term
``equivariant''.

\section{Chambers for the metric and perturbation}
\label{sec:chambers}

Let $Y$ be a $\spc$ rational homology $3$--sphere
with Riemannian metric $g$ and perturbation $1$--form $\nu$.
Let $\cd$ be the corresponding Chern--Simons--Dirac functional (which depends
on the choice of a reference spin connection over $Y$). The critical points
of $\cd$ are the solutions $(B,\Psi)$ to the following $3$--dimensional monopole
equations:
\be{align*}
*(\hatF_B+i\,d\nu)&=\frac12\si(\Psi,\Psi),\\
\prtl_B\Psi&=0.
\end{align*}
Here $B$ is a spin connection over $Y$ and $\hatF_B$ half the curvature of the
induced $\U1$--connection $\check B$, whereas $\Psi$ is a section
of the spin bundle over $Y$. For more details, including the definition of
the quadratic term $\si$, see \cite[Subsection~3.2]{Fr13}.

Any representative of the reducible critical point $\theta$
has the form $(B-i\nu,0)$ with $\check B$ flat.
We will say $(g,\nu)$ is a {\em TND-pair} if
$\theta$ is non-degenerate (i.e.\ if the kernel of the Dirac operator
$\prtl_{B-i\nu}$ is zero) and an {\em ND-pair} if all critical points of $\cd$
are non-degenerate. According to \cite[Proposition~8.1.1]{Fr13}, if $g$ is
any metric on $Y$ then $(g,\nu)$ will be an ND-pair for generic $\nu$.

To any TND-pair $(g,\nu)$
we associate an index $I(g,\nu)\in\mm(Y)$ as follows.
Let $X$ be any $\spc$ Riemannian $4$--manifold with one tubular end $\rpy$.
Let $A$ be a spin connection over $X$ whose restriction to the end agrees with
the spin connection induced by $B-i\nu$. Now set
\be{equation}\label{eqn:igv}
I(g,\nu)=\ind_{\co}(D_A)-\frac18(c_1(\cl_X)^2-\si(X)),
\end{equation}
where $D_A:L^2_1(X;\bs^+)\to L^2(X;\bs^-)$, cf.\ \cite[Chapter~9]{Fr13}.
Here $\bs^\pm$ denote the spin bundles over $X$.
Then $m:=I(g,\nu)$ depends only on $g,\nu$ and
the orientation and $\spc$ structure on $Y$, essentially because the 
right-hand side of \Ref{eqn:igv} vanishes when $X$ is closed.
In situations where we are dealing with several pairs $(g,\nu)$ simultaneously
we will often write $\theta_m$ and $\fl(Y,m)$ instead of $\theta$ and $\fl_Y$,
resp.

The following theorem is crucial to this whole paper:

\be{thm}\label{thm:realize}
For any $m\in\mm(Y)$ there exists a metric $g$ and $1$--form
$\nu$ on $Y$ such that $m=I(g,\nu)$.
\end{thm}

The idea of the proof, which occupies the remainder of this section,
is to graft onto $Y$ a family of Berger metrics
on $S^3$ and then apply Hitchin's computation \cite{Hitchin1} of the 
spectrum of the Dirac operator for these metrics together with 
B\"ar's theorem \cite{Baer2} on spectra of twisted Dirac operators
over connected sums.

Let $\mm'(Y)$ be the set of all $m\in\mm(Y)$
that can be realized as $I(g,\nu)$
for some pair $(g,\nu)$. Theorem~\ref{thm:realize} is an immediate consequence
of the following lemma.

\be{lemma}\label{lemma:mm-unbounded}
\be{enumerate}
\item[(i)]The set $\mm'(Y)$ is unbounded both below and above.
\item[(ii)]If $m_1,m_3\in\mm'(Y)$ and $m_2\in\mm(Y)$ satisfy
$m_1<m_2<m_3$ then $m_2\in\mm'(Y)$.
\end{enumerate}
\end{lemma}

The proof will involve the notion of spectral flow, which for the present 
purposes we will interpret as an index (cf.\ \cite{D5,Robbin-Salamon}).
To explain this, suppose $Z$ is a closed manifold and $E\to Z$ an Hermitian
vector bundle. Let
$P_t:\Ga(E)\to\Ga(E)$, $a\le t\le b$ be a smooth family of
first-order, self-adjoint, elliptic operators such that
the kernels of $P_a$ and $P_b$ are zero.
Choose a smooth function
$f:\R\to[a,b]$ such that $f(t)=a$ for $t\le-1$ and $f(t)=b$ for
$t\ge1$. We then define the spectral
flow $\text{SF}$ of the family $\{P_t\}_{a\le t\le b}$ to be the index of
\[\frac\prtl{\prtl t}+P_{f(t)}:L^2_1\to L^2,\]
which is independent of $f$. (It is easy to see that
$\frac\prtl{\prtl t}+P_{f(t)}$ is an elliptic operator.)
We will make use of the following simple fact:

\be{obs}\label{obs:nbna}
If $h:[a,b]\to[0,\infty)$ is a continuous function such that
$\ker(P_t-h(t))=0$ for all $t$ then
$\text{SF}=n_b-n_a$, where $n_t$ is the number of eigenvalues of $P_t$ in the
interval $(0,h(t))$, counted with multiplicity.\square
\end{obs}


{\em Proof of Lemma~\ref{lemma:mm-unbounded}:}
(i) Hitchin \cite{Hitchin1} computed the spectrum of the Dirac
operator over $S^3$ for a family of metrics $\{g_T\}_{T>0}$ known as
Berger metrics, $g_1$ being the standard metric. Let $D_T$ 
denote the Dirac operator over $S^3$ for the metric $g_T$. The eigenvalues of
$D_T$ are given by a countable family of smooth functions of $T$,
each function
describing an eigenvalue of a fixed multiplicity. Each function is either
$\ge\sqrt 2$ for all $T$, or has positive derivative and one zero. Moreover,
there are infinitely many functions of the second kind.
The set of those $T$ for which $D_T$ has non-trivial kernel is a closed
and discrete subset of $\R$ whose smallest element is $4$.

Fix $\tau>1$ with $\ker(D_\tau)=0$. The idea now is to graft the family
$\{g_T\}_{1\le T\le\tau}$ onto $Y$ and apply B\"ar's gluing theorem for
spectra of twisted Dirac operators over connected sums \cite{Baer2}. To this
end we fix $x_0\in S^3$ and construct a new family of metrics $g_{T,r}$
on $S^3$ by flattening $g_T$ in the geodesic ball of radius $r$ around $x_0$.
To make this precise, choose $r_0>0$ such that
$S^3$ contains a $g_T$--geodesic ball
$\utr$ of radius $2r_0$ around $x_0$ for $1\le T\le\tau$. 
For the remainder of this paragraph
we fix $T\in[1,\tau]$ and use $g_T$ to define exponential maps etc.
Let $\bar g_T$ denote the flat metric on $\utr$ which under the
exponential map corresponds to the
constant metric $g_T(x_0)$ on $T_{x_0}S^3$. Choose
a smooth function $\beta:\R\to\R$ such that $\beta(t)=0$ for $t\le1$ and
$\beta(t)=1$ for $t\ge2$.
For $0<r<r_0$ define $\beta_r:\utr\to\R$ by
\[\beta_r(\exp(v))=\beta(|v|/r)\]
for $v\in T_{x_0}S^3$ of norm $<2r_0$.
Let $\gtr$ be the metric on $S^3$
which agrees with $g_T$ outside $\utr$ and such that
\[\gtr=\beta_r g_T+(1-\beta_r)\bar g_T\quad\text{on $\utr$.}\]
Then $\gtr=\bar g_T$ on the ball of radius $r$ around $x_0$.

We will need to compare Dirac operators for different metrics on $S^3$ .
In general, if $g,\ti g$ are Riemannian metrics on a spin manifold $Z$, one
can identify the spin bundles such that $\ti g$--Clifford multiplication
with a tangent vector $e$ corresponds to $g$--Clifford multiplication
with $v(e)$, where $v$ is the positive $g$--symmetric bundle automorphism
of the tangent bundle $TZ$ with $v^*g=\ti g$.
Let $\nabla,\ti\nabla$ be the Riemannian connections
and $D,\ti D$ the Dirac operators on $Z$ for the metrics $g,\ti g$, resp.
Then for any spinor $s$ one has a pointwise estimate
\[|(\ti D-D)s|\le
C\left(\al\cdot|\nabla s|+(1+\al)\cdot|\ti\nabla-\nabla|\cdot|s|\right),\]
where $\al=|\ti g-g|$, norms on tensors are taken with respect to $g$, and the
constant $C$ depends only on the dimension of $Z$.

Returning to our main discussion, let $\dtr$ be the Dirac operator on
$S^3$ for the metric $\gtr$. It is a simple exercise to show that
\[\sup_{1\le T\le\tau}\|\gtr-g_T\|_{C^1}\to0\quad\text{as $r\to0^+$.}\]
Recalling the explicit formula for the Riemannian connection in terms of
the metric, we conclude that
\be{equation}\label{eqn:dtr-dt}
\sup_{1\le T\le\tau}\|\dtr-D_T\|\to0\quad\text{as $r\to0^+$,}
\end{equation}
where we use the operator norm for bounded operators $L^2_1\to L^2$.

We now introduce a suitable family of metrics on the
cylinder $\R\times S^2$. As in \cite[p\ 230]{Baer2} choose a smooth function
$\rho:\R\to\R$ such that
\be{itemize}
\item $\rho(t)=|t|$ for $|t|\ge1$,
\item $0<\rho(t)\le1$ for $|t|\le1$,
\item $|\rho'(t)|\le1$ for all $t$.
\end{itemize}
For $u>0$ set $\rho_u(t)=u\rho(t/u)$. Let $\gcyl$ denote the warped product
metric on $\R\times S^2$ given by
\[ds^2=dt^2+\rho_u(t)^2d\si^2,\]
where $d\si^2$ is the standard metric on $S^2$ and $t$ is the $\R$--coordinate.

Choose a metric $g$ on $Y$ for which $Y$
contains a Euclidean ball of radius $r_0$ around some point $y_0$.
Fix $1\le T\le\tau$ and $0<u<r<r_0$. In this
paragraph, $S^3$ and $Z_r:=(-r,r)\times S^2$ will have the metrics
$\gtr$ and $\gcyl$, respectively, and $B_u(p)$ will denote the open geodesic
ball of radius $u$ around a point $p$. Let $A^\pm$ denote the region in
$Z_r$ defined by the inequality $u<\pm t<r$. Choose an orientation
preserving isometry $\iota^+$ between $A^+$ and the corresponding annulus
in $Y$ around $y_0$. Choose a smooth family $\{\phi_T\}_{T>0}$ of
orientation reversing orthogonal transformations $\R^3\to T_{x_0}S^3$,
where $\R^3$ has the Euclidean metric. Identifying $A^-$ with an annulus in
$\R^3$ in the obvious way we obtain from $\phi_T$ an orientation preserving
isometry $\iota^-$ between $A^-$ and the corresponding annulus in $S^3$
around $x_0$. We can now introduce the glued smooth manifold
\[Y_T:=S^3\setmin\oline{B_u(x_0)}\:\cup_{\iota^-}\:(-r,r)\times S^3
\:\cup_{\iota^+}\:Y\setmin\oline{B_u(y_0)},\]
which inherits a metric $\tgtru$ from $\gtr,\gcyl,\gy$.
Set $Y^-:=Y\setmin B_{r_0}(y_0)$.

Note that $Y^-$ can also be regarded as a subset of $Y_T$.
Since $g_1$ is the standard metric on $S^3$, there is clearly a diffeomorphism
$q:Y\to Y_1$ which is the identity on $Y^-$.

The disjoint union
\[\ti Y:=\bigcup_{1\le T\le\tau}Y_T\times\{T\}\]
has a natural
smooth structure
such that the projection $\pi:\ti Y\to[1,\tau]$ is a submersion. Let $\Xi$ be a
vector field on $\ti Y$ whose restriction to $Y^-\times[1,\tau]$ is
$\frac\prtl{\prtl T}$ and such that $\pi_*\Xi=\frac\prtl{\prtl T}$.
The flow of $\Xi$
together with the diffeomorphism $q:Y\to Y_1$ yields a
diffeomorphism $Y\times[1,\tau]\to\ti Y$ which commutes with the projections
onto $[1,\tau]$ and is the identity on $Y^-\times[1,\tau]$. Let $\gtru$
be the pull-back of the metric $\tgtru$ by the induced diffeomorphism
$Y\to Y_T$.

Choose a connection in the determinant
line bundle $\cl_0\to Y$ which is flat outside $Y^-$ and such that
the kernel of the corresponding Dirac operator $D$ for the metric $g$ is zero
(see \cite[Lemma~8.1.2]{Fr13}).
In this proof, all Dirac operators over $Y$ will be defined using this
connection. Let $\dtru$ denote the Dirac operator on $Y$ for the metric
$\gtru$. Choose $\La>0$ such that
neither of the operators $D_1$ and $D$ has any 
eigenvalue in $[-\La,\La]$. Because of \Ref{eqn:dtr-dt} we can find an
$r_1\in(0,r_0)$ such that for $0<r<r_1$ the operator $D_{1,r}$ has no
eigenvalue in $[-\La,\La]$.
By the claim on \cite[p 230]{Baer2} and the explicit expression for $\delta$
given after that claim, there is a positive, increasing function $f$
with the following significance: Set
\[K=\left\{(r,u)\in\R^2\st 0<r<r_1,\;0<u<f(r)\right\}.\]
Then for $(r,u)\in K$ the kernel of $D_{1,r,u}$ is zero.
Since any two points in $K$ can be joined by a smooth path in $K$,
it follows that the spectral flow of the Dirac operator on $Y$ along a path
of metrics from $g$ to $g_{1,r,u}$ is independent of $(r,u)\in K$.

Let $T_1<\dots<T_N$ be the values of $T\in(1,\tau)$
for which $\ker(D_T)\neq0$. Set $T_0=1$, $T_{N+1}=\tau$,
$d_k=\dim \ker(D_{T_k})$.
Choose $c,\lla$ such that
\[0<c<\min_{0\le k\le N}(T_{k+1}-T_k),\qquad0<\lla\le\La,\]
and
\be{enumerate}
\item[(a)]$\ker(D_T-\lla)=0$ if $|T-T_k|\le c$ for some $k$,
\item[(b)]$D_{T_k-c}$ has no eigenvalue in $(0,\lla)$,
\item[(c)]$D_{T_k+c}$ has exactly $d_k$ eigenvalues in $(0,\lla)$,
counted with multiplicity.
\end{enumerate}
Of course we also have
\be{enumerate}
\item[(d)]$\ker(D_T)=0$ if $1\le T\le\tau$ and $|T-T_k|\ge c$ for all $k$.
\end{enumerate}
By \Ref{eqn:dtr-dt} there is an $r_2\in(0,r_1)$ 
such that if $0<r<r_2$ then (a)--(d) also hold with $D_{\bullet,r}$
in place of $D_\bullet$. Fix such an $r$. Another application
of B\"ar's gluing theorem then
shows that if $u>0$ is sufficiently small then (a)--(d) also hold with
$D_{\bullet,r,u}$ in place of $D_\bullet$.
For such $r,u$ the spectral flow of the path of operators
$\{\dtru\}_{1\le T\le\tau}$ is therefore $\sum d_k$. (This follows from
Observation~\ref{obs:nbna} and the addition formula for the index 
\cite[Corollary~C.0.1]{Fr13}.) Since this sum can be made
arbitrarily large by suitable choice of $\tau$, we conclude that
the set $\mm'(Y)$ is unbounded above. 
Now recall that reversing the orientation of
an odd-dimensional base-manifold reverses the sign of the
Dirac operator, see \cite[Section~3.1]{Fr13}.
By replacing $g_T$ with $\iota^*(g_T)$ for some fixed orientation reversing
diffeomorphism $\iota$ of $S^3$ one can therefore show in a similar way that
$\mm'(Y)$ is unbounded below.

(ii) If $(g_0,\nu_0)$, $(g_1,\nu_1)$ are any TND-pairs on $Y$ then these can
be connected by a smooth path $\{(g_t,\nu_t)\}_{0\le t\le1}$ of metrics
and $1$--forms. By essentially
the same argument as in \cite[Lemma~8.1.2]{Fr13} one can show that by replacing
$\nu_t$ by $\nu_t+t(1-t)\om$ for a generic $1$--form $\om$ on $Y$ one can
arrange that
\be{equation}\label{dimkernel}
\dim_{\co}\,\ker(D_{(g_t,\nu_t)})\le1\qquad\text{for $0\le t\le1$},
\end{equation}
where $D_{(g_t,\nu_t)}$ denotes the Dirac operator on $Y$ for the metric
$g_t$ and spin connection $B-i\nu_t$.
Using Observation~\ref{obs:nbna} it is easy to see
that this implies statement~(ii) of the lemma.
\square

\section{Moduli spaces and index}\label{sec:modsp-ind}

As in \cite[Section~1.4]{Fr13} let $X$ be a $\spc$ Riemannian $4$--manifold
with tubular ends $\rpy_j$, $j=1,\dots,r$,
where each $Y_j$ is a rational homology sphere.
Let $(g_j,\nu_j)$ be an ND-pair for $Y_j$ as defined at the beginning of
Section~\ref{sec:chambers}. By a {\em GM-pair} for $X$
(with respect to the above structure) we mean a pair $(\bar g,\mu)$ consisting
of a metric $\bar g$ and $2$--form $\mu$ on $X$ which restrict to
$1\times g_j$ resp.\ $\pi^*_jd\nu_j$ over the $j$'th end,
where $\pi_j:\rpy_j\to Y_j$.
We will refer to $I(g_j,\nu_j)$ as the {\em index} of $(\bar g,\mu)$
over the end $\rpy_j$.

Now suppose we are given such $\bar g,\mu$ as well as vectors
$\vec\al=(\al_1,\dots,\al_r)$
and $\vec\prt=(\prt_1,\dots,\prt_r)$, where $\al_j\in\fl_{Y_j}$ and
$\prt_j$ is a perturbation parameter for the Seiberg--Witten equations over
$\ry$. We define the moduli space $M(X;\vec\al)=M(X;\vec\al;\vec\prt)$
to be the space of gauge
equivalence classes of solutions $(A,\Phi)$ to the Seiberg--Witten equations
\be{equation*}
\be{gathered}
\left(\hatF_A+i\mu+\text{Pert}(S,\vec\prt)\right)^+=Q(\Phi),\\
D_A\Phi=0
\end{gathered}
\end{equation*}
such that $(A,\Phi)$ is asymptotic to $\al_j$ over the end $\rpy_j$.
Here $A$ is a spin connection over $X$ and $\hatF_A$ half the curvature of
the induced $\U1$--connection, $\Phi$ is a section of the positive
spin bundle over $X$ and $Q(\Phi)$ a certain quadratic function of $\Phi$,
whereas
$\text{Pert}(S,\vec\prt)$ is a certain compactly supported perturbation,
see \cite[Section~3.4]{Fr13} for more details.

In the case $X=\ry$ we denote by
$M(\al,\beta)$ the moduli space over $\ry$ with asymptotic limits
$\al$ at $-\infty$ and $\beta$ at $\infty$.

By definition, $\dim M(X;\vec\al)$ is the {\em expected dimension} of 
the moduli space, i.e.\ the index of any Fredholm operator of a certain kind,
see \cite[p.\ 33]{Fr13}.
In the case of reducible
limits $\vec\theta=(\theta,\dots,\theta)$ that Fredholm operator can be
taken to be the one induced by
$(-d^*+d^+)\oplus D_A$, which yields
\be{equation}\label{eqn:dimMred}
\dim M(X;\vec\theta)=-b_0(X)+b_1(X)-b^+(X)+2\,\ind_{\co}(D_A).
\end{equation}
Here $\ind_{\co}(D_A)$ is as in \Ref{eqn:igv}. (A priori, one should consider
$D_A$ as a map between certain $L^p$ rather than $L^2$ Sobolev-spaces, but this
gives the same index.)

The addition formula for the index (see \cite[Corollary~C.0.1]{Fr13})
tells us how the expected dimension
of the moduli space changes when two ends of $X$ with the same asymptotic limit
are being glued.
In the simplest case when $X$ has only two ends $\R_\pm\times Y$ one has
\be{equation}\label{eqn:addnform}
\dim M(\xt)=\dim M(X;\al,\al)+n_\al,
\end{equation}
where $\xt$ is the manifold obtained from $X$ by replacing the two ends
by a neck $[-T,T]\times Y$, and $n_\al=0$ if $\al$ is irreducible and
$n_\al=1$ otherwise. An analogous formula holds if $X$ has further
ends which are not being glued.

Now let $Y$ be any $\spc$ rational homology $3$--sphere equipped with
an ND-pair of index $m$. We will associate an index (or degree)
to any $\al\in\fl_Y$.
Let $X$ be any connected $\spc$ Riemannian $4$--manifold with
one tubular end $\rpy$. We define
$\ind(\al)\in2\mm(Y)+\z$ and $\ind_2(\al)\in\z/2\z$ by
\be{equation}\label{dim-formula}
\be{gathered}
\dim M(X,\al)=
\ind(\al)+\frac14\left(c_1(\cl_X)^2-\si(X)\right)-n_\al+b_1(X)-b^+(X),\\
\dim M(X,\al)\equiv
\ind_2(\al)-n_\al+b_1(X)-b^+(X)\mod2.
\end{gathered}
\end{equation}
Thus, $\ind_2(\al)\equiv\ind(\al)\mod2$ if $Y$ is an integral homology sphere.
Moreover, by \Ref{eqn:dimMred} one has
\[\ind(\theta)=2m,\qquad\ind_2(\theta)=0.\]
If $\al$ is irreducible then an application of the addition formula
\Ref{eqn:addnform} yields
\be{gather*}
\ind_2(\al)\equiv\dim M(\theta,\al)\mod2.
\end{gather*}

The reason for our preference for Floer cohomology rather than homology
is that we want $\ind(\al)$ to appear with a positive sign in
\Ref{dim-formula}, in which case the Seiberg--Witten invariants of a 
$4$--manifold with boundary lie in the Floer {\em co}homology of the
boundary.

\section{Perturbations and reducibles}
\label{sec:pert-red}

We will usually assume that
$\bar g,\mu,\vec\prt$ have a certain generic form with regard to reducibles
in $M(X;\vec\theta)$. To explain this, let $\bar g$ be given.
We first take $\vec\prt=0$. If $b^+(X)>0$
choose $\mu$ such that $M(X;\vec\theta)$ is free of reducibles.
If $b_1(X)=0=b^+(X)$ then $M(X;\vec\theta)$ contains a unique 
reducible point for any $\mu$, and if in addition
$\dim M(X;\vec\theta)\ge-1$ then by \cite[Lemma~14.2.1]{Fr13}
we may choose $\mu$
such that the reducible point is regular.
Having chosen $\bar g,\mu$ we will now show that these properties of 
the moduli space are retained if we allow small perturbation parameters
$\prt_j$.

Let $\Prt_j$ be a Banach space of perturbation parameters $\prt_j$ as in 
\cite[Section~8.2]{Fr13}, and set $\Prt=\Prt_1\times\dots\Prt_r$.

\be{prop}\label{prop:pert-red}
In the above situation the following hold when
$\|\vec\prt\|$ is sufficiently small:
\be{enumerate}
\item[(i)]If $b^+(X)>0$ then $M(X;\vec\theta)$ is free of reducibles.
\item[(ii)]If $b_1(X)=0=b^+(X)$
then $M(X;\vec\theta)$ contains a unique reducible point $\om(\vec\prt)$.
Moreover, if $\dim M(X;\vec\theta)\ge-1$
then this point is a regular point of the moduli space.
\end{enumerate}
\end{prop}

\proof We prove (ii) and leave the easier part (i) to the reader. We can
identify the reducible part of $\cb(X;\vec\theta)$ with a configuration space
$\cb(\cl_X)$ of connections in $\cl_X$. Let $A_0$ be a connection in $\cl_X$
representing $\om(0)$. By \cite[Proposition~2.3.1]{Fr13} any point in 
$\cb(\cl_X)$ has a representative $A$ satisfying $d^*a=0$, where $a=A-A_0$.
Moreover, $[A]$ lies in the reducible part $M_{\vec\prt}\subset\cb(\cl_X)$
of the moduli space $M(X;\vec\theta;\vec\prt)$ if and only if
\[0=\frac12F^+_A+i\mu^++\eps(A,\vec\prt)
=\frac12d^+a+\eps(A,\vec\prt),\]
where $\eps$ is a perturbation term. For $\del>0$ set
\[U_\del=\{a\in\llw pw1(X;i\La^1)\st d^*a=0,\;\|a\|_{\llw pw1}<\del\},\]
and consider the smooth map
\[f:U_1\times\Prt\to\lw pw,\quad(a,\vec\prt)\mapsto
\frac12d^+a+\eps(A_0+a,\vec\prt).\]
Since $b_1(X)=0=b^+(X)$, the derivative in the first variable $D_1f(0,0)=d^+$
is an isomorphism. By the implicit function theorem there exist $\eps,\del>0$
such that if $\|\vec\prt\|<\eps$ then $f$ has a unique zero
$(a(\vec\prt),\vec\prt)$ with $a(\vec\prt)\in U_\del$. We will show that
$[A_0+a(\vec\prt)]$ is the only point in $M_{\vec\prt}$ when $\vec\prt$
is sufficiently small. So suppose $\{\vec\prt(n)\}$ is a sequence in $\Prt$
converging to $0$, and $\om_n\in M_{\vec\prt(n)}$ for each $n$.
After passing to a
subsequence we may assume that $\om_n$ chain-converges. (This does not quite
follow from \cite[Theorem~1.3.1]{Fr13} since we vary the
perturbation parameters,
but the proof of that theorem works just as well in the present situation.)
For a rational homology $3$--sphere there is only one reducible point in
$\ti\fl=\fl$, so the chain-limit lies in $M_0$ and is therefore equal
to $[A_0]$. By \cite[Proposition~6.4.1]{Fr13} the sequence $\{\om_n\}$ actually
converges in $\cb(\cl_X)$. Now, $A_0+U_1\to\cb(\cl_X)$ is an open map, so
when $n$ is sufficiently large there is a representative of $\om_n$ of the form
$A_0+a_n$ such that $d^*a_n=0$ and $\|a_n\|_{\llw pw1}\to0$. But then
$a_n=a(\vec\prt(n))$ for $n\gg0$.\square

\section{Orientations}

Let $Y$ be a $\spc$ rational homology $3$--sphere equipped with an
ND-pair. By an {\em orientation} of a critical point $\al\in\fl^*_Y$
we mean an orientation of $\cb(\theta,\al)=\cb(\ry;\theta,\al)$,
or more precisely, a section of the orientation cover $\lla\to\cb(\theta,\al)$
constructed in \cite[Section~4.2]{Fr13}.

We will now state the conventions that we will use for orienting
moduli spaces. Let $X$ be a homology oriented Riemannian
$\spc$ $4$--manifold with tubular ends modelled on rational homology spheres. 
The {\em canonical} orientation $\ocan$
of $\cb(X;\vec\theta)$ is the one determined by the homology orientation
of $X$ as in \cite[Proposition~12.2.1]{Fr13}. 
Now suppose the ends of $X$ are $\R_-\times Y^-_j=\R_+\times(-Y^-_j)$,
$j=1,\dots,r^-$ and
$\R_+\times Y^+_j$, $j=1,\dots,r^+$. We will refer to these as ``negative''
and ``positive'' ends, respectively.
Let $\al^\pm_j$ be a monopole over $Y^\pm_j$
which is irreducible and equipped with an orientation $o^\pm_j$ for
$j=1,\dots,\rho^\pm$ and reducible for $j=\rho^\pm+1,\dots,r^\pm$.
Let $\vec\al^\pm:=(\al^\pm_1,\dots,\al^\pm_{r^\pm})$. We give
$\cb(X;\vec\al^-,\vec\al^+)$ the orientation $o$ satisfying
\[o^-_1\#\cdots\#o^-_{\rho^-}\#o=\ocan\#o^+_1\#\cdots\#o^+_{\rho^+}\]
as orientations of $\cb(X;\vec\theta,\vec\al^+)$, where $\#$ is as defined
in \cite[Section~13.6]{Fr13}.
Note that this definition is sensitive to the orderings of the
(positive and negative) ends of $X$. We call $o$ the
{\em proper} orientation of $\cb(X;\vec\al^-,\vec\al^+)$. Recall that
an orientation of $\cb(X;\vec\al^-,\vec\al^+)$ determines an orientation
of the regular and irreducible part of the moduli space
$M(X;\vec\al^-,\vec\al^+)$.

If $X=\ry$, where the negative and positive ends are the obvious ones,
then we obtain orientations of $\cb(\ry;\al,\theta)$ and
$\cb(\ry;\al,\beta)$ for all oriented $\al,\beta\in\fl^*_Y$.
In particular, an orientation of $\al\in\fl^*_Y$ determines an orientation
of the corresponding element of $\fl^*_{-Y}$
(which in turn induces the original orientation of $\al$, essentially 
because either $M(\al,\theta)$ or $M(\theta,\al)$ is even-dimensional).

Returning to the general case, note that changing
the status of an ``irreducible''
end of $X$ from negative to positive alters the
orientation of $\cb(X;\vec\al^-,\vec\al^+)$ by a sign which only depends on
the parity of $b_0(X)+b_1(X)+b^+(X)$, the mod~$2$ indices
$\ind_2(\al^\pm_j)$, and the new and old orderings of the ends of $X$.

\be{prop}\label{prop:glue-ortn1}
Let $X_1$ and $X_2$ be two $4$--manifolds as above, both connected.
Suppose $X_1$ has
only one positive end $\R_+\times Y$, with limit $\al$, and that $X_2$
has only one negative end $\R_-\times Y$, with the same limit $\al$.
Let $\vec\al^\pm$ be the limits assigned to the negative ends of $X_1$
and the positive ends of $X_2$, respectively. For $T>0$ let $\xt$ denote
the glued manifold obtained from the disjoint union
$X_1\cup X_2$ by replacing the two
ends $\R_\pm\times Y$ with a neck $[-T,T]\times Y$. We regard $\xt$
as having the same negative (resp.\ positive) ends as $X_1$ (resp.\ $X_2$).
Let $\xt$ have
the glued homology orientation defined in \cite[Section~13.6]{Fr13}.
Let each limit $\al,\al^\pm_j$ be oriented if irreducible,
and let $o_1,o_2,o'$ denote the proper orientations
of $\cb(X_1;\vec\al^-,\al)$, $\cb(X_2;\al,\vec\al^+)$, and
$\cb(\xt;\vec\al^-,\vec\al^+)$, respectively. Then
\be{equation}\label{eqn:o1o2o}
o_1\#o_2=o'.
\end{equation}
If $\al$ is irreducible, then any ungluing map as in
\cite[Theorem~11.1.1]{Fr13} is an orientation
preserving diffeomorphism from an open subset of $M(\xt;\vec\al^-,\vec\al^+)$
onto an open subset of $M(X_1;\vec\al^-,\al)\times M(X_2;\al,\vec\al^+)$.
\end{prop}

\proof Equation~\ref{eqn:o1o2o} is easily deduced from the properties
of the $\#$--operation mentioned in \cite[Section~13.6]{Fr13}. The statement
about the ungluing maps then follows from \cite[Theorem~13.4.1]{Fr13} by
unraveling the definition of $\#$.
\square

We call a diffeomorphism {\em $\eps$--preserving} (where $\eps=\pm1$)
if it preserves or reverses
orientations according to the sign of $\eps$.

\be{prop}\label{prop:glue-ortn2}
Let $X$ be a $4$--manifold as above which is connected and has
only one negative end
$\rmy$ and only one positive end $\rpy$. Let the same oriented irreducible
limit $\al$ be assigned to both ends. For $T>0$ let $\xt$ denote
the glued manifold obtained from 
$X$ by replacing the two ends $\R_\pm\times Y$ with a neck $[-T,T]\times Y$.
Let $\xt$ have the glued homology orientation defined in
\cite[Section~13.5]{Fr13}. Let $\cb(X;\al,\al)$ have the proper orientation,
and let $\ti o$ denote the corresponding glued orientation of $\cb(\xt)$
defined in \cite[Section~13.4]{Fr13}. Then 
\[\ti o=(-1)^{ab+a+b}\ocan,\]
where $a=\ind_2(\al)$ and $b=\dim M(\xt)$.
Any ungluing map as in \cite[Theorem~11.1.1]{Fr13} is a
$(-1)^{ab+a+b}$--preserving diffeomorphism from an open subset of $M(\xt)$
onto an open subset of $M(X;\al,\al)$.
\end{prop}

\proof The sign requires a simple computation, which we omit. The 
statement about the ungluing map is just a special case of
\cite[Theorem~13.4.1]{Fr13}.\square

\section{Irreducible Floer groups}
\label{sec:irrflg}

Let $Y$ be a $\spc$ rational homology $3$--sphere and $m\in\mm(Y)$. We will
now define the irreducible Floer cohomology group $\hf^*(Y,m)$.
By Theorem~\ref{thm:realize} and \cite[Proposition~8.1.1]{Fr13}
we can find an ND-pair $(g,\nu)$ on $Y$ of index $m$.

We first explain how to associate signs to index~$1$ flow-lines.
Let $\al,\beta\in\fl:=\fl_{(g,\nu)}$, not both reducible, and oriented if 
irreducible. Let $\R$ act on
$M:=M(\al,\beta)$ by
\be{equation}\label{eqn:rwtr}
r\cdot\om:=\ct^*_{-r}\om,
\end{equation}
for $(r,\om)\in\R\times M$, where the translation operator $\ct_s$
on $\ry$ is defined by $\ct_s(t,y)=(t+s,y)$. If
$\dim M=1$ then the elements of the orbit space
\[\chm:=M/\R\]
are precisely the connected components of $M$,  and to each such component $\Ga$
we associate a sign $\pm1$ according as to whether the diffeomorphism
$\R\to\Ga$, $r\mapsto r\cdot\om$ preserves or reverses orientation for
$\om\in\Ga$. 

Let $\fl^q$ be the set of all elements of $\fl^*$ of index $q$.
Let $\cf^q$ be the Abelian group generated by all pairs $(\al,o)$ where
$\al\in\fl^q$ and $o$ is an orientation of $\al$, subject to the relation
$(\al,-o)=-(\al,o)$. For convenience of notation, however, we will always
assume that an orientation has been chosen for every element of $\fl^*$,
in which case $\cf^q$ can be identified with the free Abelian group
generated by $\fl^q$.

We define the differential $d=d^q:\cf^q\to\cf^{q+1}$ as usual by 
\[d\al=\sum_\beta(\#\chm(\al,\beta))\beta,\]
where the sum is taken over all $\beta\in\fl^*$ of index $q+1$.
Here and elsewhere in this paper, $\#$ denotes a signed count.
Note that the differential depends on the choice
of a small perturbation of the Seiberg--Witten equations over $\ry$
as in \cite{Fr13}. However, this perturbation will not be reflected
in our notation.

We also need to study flow-lines running between irreducibles and the
reducible. Noting that
\be{align*}
\dim M(\al,\theta_m)&=2m-1-\ind(\al),\\
\dim M(\theta_m,\al)&=\ind(\al)-2m,
\end{align*}
we define, following \cite{D5}, homomorphisms
\be{align*}
&\del:\cf^{2m-2}\to\z,\quad\al\mapsto\#\chm(\al,\theta_m),\\
&\del':\z\to\cf^{2m+1},\quad1\mapsto\sum_\al(\#\chm(\theta_m,\al))\al,
\end{align*}
where the sum is taken over all $\al\in\fl^*$ of index $2m+1$.

\be{prop}$d\circ d=0$, $\del\circ d=0$, $d\circ\del'=0$.
\end{prop}
\proof This follows from \cite[Section~12.2]{Fr13} and
Proposition~\ref{prop:glue-ortn1}.\square

For a fixed Abelian group $G$ let $\hf^*$ denote the cohomology of
$\cf^*\otimes G$, and let
\be{equation}\label{eqn:deldeldef}
\hf^{2m-2}\oset{\del_0}\longrightarrow G\oset{\del'_0}
\longrightarrow\hf^{2m+1}
\end{equation}
be the homomorphisms induced by $\del,\del'$. It is convenient to
extend $\del_0$ to all of $\hf^*$ by setting $\del_0:=0$ in degrees
different from $2m-2$, and similarly for $\del$.

We wish to understand the dependence of the
group $\hf^*$ on $(g,\nu)$.
(As we will see later, this problem has very much to do with $\del_0,\del'_0$.)
So let $(g_j,\nu_j)$, $j=1,2$ be two ND-pairs for $Y$
and let $m_j,\fl^*_j,\cf_j,\hf_j$ be the corresponding data as introduced
above. We will investigate for which $q$ the usual cobordism construction
yields a homomorphism $\hf_2^q\to\hf_1^q$.
Let $\bar g$ be a Riemannian metric and $\mu$ a two-form on $\ry$
such that $(\bar g,\mu)$ agrees with $(1\times g_2,\pi_1^*(d\nu_2))$ 
in the region $(-\infty,-1]\times Y$ and agrees with 
$(1\times g_1,\pi_2^*(d\nu_1))$ in $[1,\infty)\times Y$. Let $X_{21}$ denote
the manifold $\ry$ equipped with the pair $\bar g,\mu$.

\be{lemma}
Let $\al\in\fl_2^*$ and $\beta\in\fl^*_1$ be critical points of
index $q$. If either $q\le2m_1$ or $q\ge2m_2-1$ then $M(X_{21};\al,\beta)$ is
compact.
\end{lemma}
\proof Since $M(X_{21};\al,\beta)$ has dimension zero, it can only fail to be
compact if it contains a sequence which has a chain-limit involving some
monopole of negative index. The only candidate for such monopole
would be the reducible point in $M(X_{21};\thet2,\thet1)$.
For such a chain-limit
to exist, however, the moduli spaces $M(\al,\thet2)$ and $M(\thet1,\beta)$
would both need to have positive dimension. The lemma now follows from
the dimension formulae for these spaces given above.\square

For any $q$ satisfying the assumptions of the previous lemma we define
the homomorphism
\[k^q:\cf^q_2\to\cf^q_1,\quad \al\mapsto\sum_\beta(\#M(X_{21};\al,\beta))\beta.\]

\be{lemma}\label{lemma:KqHF}
If $q\le 2m_1-1$ or $q\ge2m_2$ then $k^{p+1}d^p=d^pk^p$
for $p=q-1,q$, hence $k^q$ induces a homomorphism
\[K^q:\hf_2^q\to\hf_1^q.\]
\end{lemma}

\proof Again the point is that the assumptions on $q$ prevent factorizations
through the reducible monopole on $X_{21}$.\square

Now suppose $m_1\le m_2$. Reversing the roles of $Y_1$ and $Y_2$ 
we obtain a cobordism $X_{12}$ which by the previous lemma induces
a homomorphism
\[J^q:\hf^q_1\to\hf^q_2\]
for all $q$. (That $J^q$ is well-defined can also be seen by a dimension count.)

\be{prop}\label{K-inv}
If $q\le2m_1-1$ or $q\ge2m_2$ then $J^q$ and $K^q$ are inverse maps.
\end{prop}
\proof Gluing the positive end of $X_{12}$ with the negative end of $X_{21}$
yields a cobordism $\xt$ with a band of length $2T$ where the metric and
perturbation $2$--form agrees with $(1\times g_2,\pi_1^*(d\nu_2))$. Let
$\al,\beta$ be two elements of $\fl^*_1$ of index $q$. The 
assumptions on $q$ mean precisely that one of the moduli
spaces $M(X_{12};\al,\thet2)$ and $M(\thet1,\beta)$ have negative dimension. 
Therefore no sequence $\om_n\in M(\xtn;\al,\beta)$ with $T_n\to\infty$
can have a chain-limit involving the reducible monopole over
$X_{21}$. Consequently, no factorizations through reducible critical points
may appear in such chain-limits. Thus any chain-limit must lie in
$M(X_{12};\al,\ga)\times M(X_{21};\ga,\beta)$ for some $\ga\in\fl^*_2$ of
index $q$. The gluing theorem in \cite{Fr13} then says that for 
large $T$ the endomorphism $\ell^q$ of $\cf_1^q$ induced by
the cobordism $\xt$ is equal to the composite of the homomorphisms
$\cf^q_1\to\cf^q_2\to\cf^q_1$ induced by $X_{12}$ and $X_{21}$.
Considering $1$--parameter families of moduli spaces over
$\ry$ obtained by
deforming the metric and perturbation form on $\xt$
to $(1\times g_1,\pi_1^*(d\nu_1))$
we obtain a chain-homotopy between $\ell^q$ and the identity.
Hence $K^qJ^q=\id$, and a similar argument proves $J^qK^q=\id$.
\square  

A first corollary of this proposition is that if $m_1=m_2$ then
$\hf_1^*$ and $\hf_2^*$ are canonically isomorphic. 

\be{defn}For any $m\in\mm(Y)$ let $\hf^*(Y,m;G)$ denote the cohomology group
of the cochain complex $\cf^*\otimes G$ defined
using an ND-pair $(g,\nu)$ of index $m$.
\end{defn}

Similarly, we define the Floer homology $\hf_*(Y,m;G)$ as the homology of
the chain complex $\text{Hom}(\cf^*,G)$.

We set $\hf^*(Y,m):=\hf^*(Y,m;\z)$.
We will often write $\cf^*(Y,m)$ for $\cf^*$, although $\cf^*$ really
depends on $g,\nu$ and not just on $m$.

Proposition~\ref{K-inv} now gives:

\be{prop}\label{prop:irred-iso}
Let $m_1,m_2\in\mm(Y)$ and $m_1\le m_2$.
If either $q\le 2m_1-1$ or $q\ge2m_2$ then the natural
homomorphism $\hf^q(Y,m_1;G)\to\hf^q(Y,m_2;G)$ is an isomorphism.
\square
\end{prop}

Let $\tilde g$ be the round metric on $S^3$. Because $\ti g$ has positive
scalar curvature we have $\fl_{(\ti g,0)}=\{\theta\}$. Since
$I(\tilde g,0)=0$ (see
\cite[Ch.\ 9]{Fr13}), it follows that
\[\hf^*(S^3,0;G)=0.\]
Combining this with Proposition~\ref{prop:irred-iso} yields:

\be{prop}\label{s3vanish}
If $m\ge0$ and either $q\le-1$ or $q\ge2m$ then $\hf^q(S^3,m;G)=0$.
The same conclusion holds if $m\le0$ and either $q\le2m-1$ or $q\ge0$.
\square
\end{prop}


\section{The $u$--map and cobordisms}
\label{sec:umapcob}

We begin by defining, in analogy with \cite[Section~3]{Fr3},
a natural degree~$2$ homomorphism $u:\hf^q(Y,m;G)\to\hf^{q+2}(Y,m;G)$
for $q\neq2m-2,2m-1$.
The constraint on the degree is related to interaction with the reducible
critical point.
Set $\cf^*:=\cf^*(Y,m)$.
The $u$--map will be induced by a homomorphism
\[v:\cf^*\to\cf^{*+2}\]
defined in all degrees. Let $\al,\beta\in\fl(Y,m)$
be distinct critical points, and let $\e\to M(\al,\beta)$
and $\f\to\cb^*([-1,1]\times Y)$ be the natural complex line bundles
associated to the base-point $(0,y_0)$. (Here we use the same convention
as in \cite[Section~15.1]{Fr13}.) Let
$R:M(\al,\beta)\to\cb^*([-1,1]\times Y)$ be the natural restriction map.
This is well defined because of the unique continuation
property of Dirac operators (see \cite{JKazdan1}).
Note that we can identify $R^*\f=\e$.
Choose a generic section $s$ of $\f$ and let $\si$ be the induced section
of $\e$. If $\al,\beta$ are irreducible and
$\ind(\beta)-\ind(\al)=2$ then $\si\inv(0)$ is a finite set of
oriented points, and we define the matrix coefficient
$\la v\al,\beta\ra$ by
\[\la v\al,\beta\ra=\#\si\inv(0).\]

\be{prop}\label{dvvd}
$dv-vd+\del'\del=0$.
\end{prop}

\proof This is similar to the proof of \cite[Theorem~4(ii)]{Fr3}.
Note however, that the modifications to the sections $\si_j$ made in
part~(I) of that proof are unnecessary, as their role in ruling out
cases (b), (c) can be replaced by a transversality argument.
(See Part~(V), Case~1 of the proof of Lemma~\ref{del-k} below.)
This simplifies the proof somewhat. 
In particular, a general unique continuation
property for monopoles is not needed for the proof of the proposition.
(We do not know whether this property holds for the 
perturbed Seiberg--Witten equations used in this paper (see
\cite[Section~8.2]{Fr13}), because these equations are not of 
gradient flow type -- in fact the perturbations
are non-local in the time direction.
In particular, we do not see how to
apply the unique continuation results of \cite{KM5}.)

Up to sign, the coefficient of $\del'\del$ in the proposition is the 
Euler number of the rank~$1$ Hermitian vector bundle over
$S^2=D^2\cup_{S^1}D^2$ whose ``clutching map'' $S^1\to\U1$ has degree~$1$.
The sign can be determined using Proposition~\ref{prop:glue-ortn1}.\square

\be{defn}For $q\neq2m-2,2m-1$ we define 
\be{equation*}
u:\hf^q(Y,m;G)\to\hf^{q+2}(Y,m;G)
\end{equation*}
to be the homomorphism induced by $v$.
\end{defn}

It will follow from Proposition~\ref{wu-comm} below that $u$ is independent of
the choice of section $\si$ and of the metric and perturbations used in
the definition of the Floer cohomology.

If $W$ is any compact $\spc$ $4$--manifold with boundary $Y$ then we denote
by
\[\hw:=W\cup_Y([0,\infty)\times Y)\]
the result of attaching a half-infinite cylinder to each boundary component
of $W$. Now suppose $Y=(-Y_1)\cup Y_2$, where each $Y_j$ is a
rational homology $3$--sphere.
If $\al_j\in\fl(Y,m_j)$, $j=1,2$, then
\[\dim M(\hw;\al_1,\al_2)=\ind(\al_2)-\ind(\al_1)-n_{\al_2}+d,\]
where
\be{equation}\label{eqn:ddef}
d=\frac14(c_1(\cl_W)^2-\si(W))+b_1(W)-b^+(W).
\end{equation}
Let $m_j\in\mm(Y_j)$ and set $k=m_2-m_1+\frac d2\in\frac12\z$, so that
\[\dim M(\hw;\theta_{m_1},\theta_{m_2})=2k-1.\]
If either $b^+(W)>1$ or $k>-1$ then $W$ gives rise to a cochain map
\be{gather*}
W^\#:\cf^*(Y_1,m_1)\to\cf^{*-d}(Y_2,m_2),\\
\al\mapsto\sum_\beta(\#M(\hw;\al,\beta))\beta,
\end{gather*}
which in turn induces a homomorphism
\be{equation*}
W^*:\hf^*(Y_1,m_1;G)\to\hf^{*-d}(Y_2,m_2;G).
\end{equation*}
Note that $W^*$ only depends on the $\spc$ manifold $W$.
The constraint on $b^+$ or $k$ was imposed to rule out factorizations
through a reducible monopole on $W$. If $b^+(W)>1$ then generically there
is no reducible on $X$ even in a one-parameter family of metrics and
$2$--forms $\mu$ (cf.\ Proposition~\ref{prop:pert-red}~(i)). If $k>-1$ and
$\al,\beta$ are irreducible then for dimensional reasons no reducible monopole
over $W$ can appear in the chain-limit of a sequence in $M(W;\al,\beta)$
provided that moduli space has dimension $\le1$, and the analogous statement
holds for parametrized moduli spaces of dimension $\le1$.

These maps $W^*$ are functorial with respect to composition of cobordisms
as long as the composite cobordism has the glued homology orientation
and also satisfies $b^+>1$ or $k>-1$.

If $b^+(W)\le1$ and $k\le-1$ then $W^*$ can still be defined under certain
restrictions on the degree, cf.\ Lemma~\ref{lemma:KqHF}.

We will now show that $W^*$ commutes with the $u$--maps, in as far as
these are defined.
Let $v_j:\cf^*(Y_j,m_j)\to\cf^{*+2}(Y_j,m_j)$ be a homomorphism as above.

\be{prop}\label{wu-comm}
If $b^+(W)>1$ or $k\ge0$ then
there exists a homomorphism $\phi:\cf^*(Y_1,m_1)\to\cf^{*-d+1}(Y_2,m_2)$
such that
\[W^\#v_1-v_2W^\#+\del'\del_W+\del'_W\del=d\phi\pm\phi d\]
as maps $\cf^*(Y_1,m_1)\to \cf^{*-d+2}(Y_2,m_2)$,
where
\be{align*}
\del_W&:\cf^{2m_2-1+d}(Y_1,m_1)\to\z,\quad\al\mapsto\#M(\hw;\al,\thet2),\\
\del'_W&:\z\to\cf^{2m_1-d}(Y_2,m_2),\quad
1\mapsto\sum_\beta(\#M(\hw;\thet1,\beta))\beta,
\end{align*}
and $\del_W:=0$ in degrees different from $2m_2-1+d$.
\end{prop}

\proof See \cite[Theorem~6]{Fr3}.\square

\be{cor}\label{cor:wu-comm}
\[W^*u=uW^*:\hf^q(Y_1,m_1)\to \hf^{q-d+2}(Y_2,m_2)\]
for all $q$ for which both $u$--maps are defined.
\end{cor}

The assumption that $b^+(W)>1$ or $k\ge0$ is only used to
rule out that a $2$--dimensional moduli space $M(\hw;\al,\beta)$ with
$\al,\beta$ irreducible may contain a chain-convergent sequence whose limit
involves a reducible monopole on $\hw$. Therefore, if $b^+(X)\le1$ and
$k<0$ then the conclusion of the proposition still holds under certain
restrictions on the degree.


We will now define the $u$--maps on the equivariant Floer groups
$\uhf^*(Y;G)$ and $\ohf^*(Y;G)$ (which for brevity we will here denote
by $\uhf^*$ and $\ohf^*$). To define these maps in degree $q$ choose 
$m_1,m_2\in\mm(Y)$ with $2m_1\le q\le2m_2-3$. We define
$u:\uhf^q\to\uhf^{q+2}$ and $u:\ohf^q\to\ohf^{q+2}$ by requiring that
the left- and right-most squares in the following diagram commute:
\be{equation}\label{eqn:equmap}
\be{array}{ccccccc}
\uhf^q & \oset\approx\longrightarrow & \hf^q_{m_1} & \longrightarrow &
\hf^q_{m_2} & \oset\approx\longrightarrow & \ohf^q\\
u\downarrow\hphantom{u} & & u\downarrow\hphantom{u} & & 
u\downarrow\hphantom{u} & & u\downarrow\hphantom{u} \\
\uhf^{q+2} & \oset\approx\longrightarrow & \hf^{q+2}_{m_1} & \longrightarrow &
\hf^{q+2}_{m_2} & \oset\approx\longrightarrow & \ohf^{q+2}
\end{array}\end{equation}
Here $\hf^p_m=\hf^p(Y,m;G)$, and the horizontal maps are the natural ones.
(The maps marked as isomorphisms are so by Proposition~\ref{prop:irred-iso}).
This definition is independent of the choice of $m_1,m_2$.
Since the middle square commutes, the homomorphism $J_Y:\uhf^*\to\ohf^*$
commutes with the $u$--maps.

Let $W$ be as above, without any constraint on $b^+(W)$.
We define $\upsi(W)$ in degree $q$ by requiring that the diagram
\[\be{array}{ccc}
\hf^q(Y_1,m_1;G) & \oset{W^*}{\uset{\vphantom{\upsi(W)}}\longrightarrow} & \hf^{q-d}(Y_2,m_2;G)\\
\approx\uparrow\hphantom{\approx} & & \approx\uparrow\hphantom{\approx} \\
\uhf^q(Y_1;G) & \oset{\upsi(W)}\longrightarrow & \uhf^{q-d}(Y_2;G)
\end{array}\]
commute,
where $m_1,m_2$ are chosen so that $2m_1\le q$ and $2m_2\le q-d$. This
definition is independent of the choice of $m_1,m_2$. 
In a similar fashion we define $\opsi(W)$ and the invariants $\upsi(Z)$,
$\opsi(Z)$ of Subsection~\ref{subsec:inv4mfld} when $Z$ has one
boundary component.

We can now easily prove Theorem~\ref{thm:basic} by observing that 
when stretching $Z$ along $Y$ and using an ND-pair of index $m\gg0$ on $Y$
then no complications with reducibles can arise.

\section{Cobordisms and reducibles}

Combining the maps $\del_0,\del'_0$ of \Ref{eqn:deldeldef} with the $u$--map
we obtain, for $n\ge0$, homomorphisms
\[\del_0u^n:\hf^{2(m-n-1)}(Y,m;G)\to G,\quad
u^n\del_0':G\to\hf^{2m+2n+1}(Y,m;G)\]
which measure interaction between
reducible and irreducible critical points and will play a central role in 
this paper. We will now give an alternative description of these maps.

First some notation.
Let $(g,\nu)$ be an ND-pair for $Y$ of index $m$ and consider the corresponding
translationary invariant monopole equations on $\ry$.
Let $J$ be an open interval and $S$ a monopole over
$J\times Y$ such that $\prtl_t\csd{}(S_t)<0$,
where $S_t=S|_{\{t\}\times Y}$. (This will be the case unless $S$ is
gauge equivalent to the monopole determined by some 
critical point of $\cd$, see \cite[Section~7.1]{Fr13}.)
Let $\eps$ be a real number. If there exists a
$\tau=\tau(S)\in J$ such that $[\tau-1,\tau+1]\subset J$ and
\be{equation}\label{tau-def}
\cd(S_\tau)=\cd(\theta_m)+\eps,
\end{equation}
then we call $\tau$ the {\em $\eps$--slice} of $S$. 


Now set $K:=[-1,1]\times Y$ and let $\bl\to\cb^*(K)$ be the
natural complex line bundle.
Choose a generic smooth section $s$ of $\oplus^n\bl$. Choose $\eps>0$
such that there is no monopole $\beta$ over $Y$ with
$0<|\cd(\beta)-\cd(\theta_m)|\le\eps$. Then for any $\al\in\fl^*(Y,m)$ and
$\om\in M(\al,\theta_m)$ the $\eps$--slice $\tau(\om)$ is well-defined.
\cite[Proposition~3.4.2 and Lemma~3.3.1]{Fr13}
and an application
of the implicit function theorem show that $\tau$ is a smooth function
on $M(\al,\theta_m)$. Set
\[M'(\al,\theta_m):=\{\om\in M(\al,\theta_m)\st\tau(\om)=0,\:s(\om|_K)=0\}.\]
If $\al$ has index $2(m-n-1)$ then $M'(\al,\theta_m)$ consists of a finite
number of oriented points; let $\bar\del_n(\al)$ denote the number of these
points counted with sign. By considering the ends of the $1$--dimensional
spaces $M'(\beta,\theta_m)$ where $\ind(\beta)=2m-2n-3$ one finds that
$\bar\del_nd=0$.
Similarly, for any $\al\in\fl^*(Y,m)$ and $\om\in M(\theta_m,\al)$
the $(-\eps)$--slice $\tau(\om)$ is well-defined. Setting
\[M'(\theta_m,\al)=\{\om\in M(\theta_m,\al)\st\tau(\om)=0,\:s(\om|_K)=0\}\]
we obtain a homomorphism
\[\bar\del'_n:\z\to\cf^{2m+2n+1}(Y,m),\quad
1\mapsto\sum_\al(\#M'(\theta_m,\al))\al\]
satisfying $d\bar\del'_n=0$. Let
\[\hf^{2(m-n-1)}(Y,m;G)\oset{\del_n}\longrightarrow G\oset{\del'_n}
\longrightarrow\hf^{2m+2n+1}(Y,m;G)\]
be the homomorphisms induced by $\bar\del_n,\bar\del'_n$.
This agrees with the previous definition in the case $n=0$.

\be{prop}\label{prop:delnun}
$\del_n=\del_0u^n$ and $\del'_n=u^n\del'_0$ for all $n\ge0$.
\end{prop}
\proof Moving one base-point towards $-\infty$ (see the proof
of Theorem~6 in \cite{Fr3}) one finds
that $\del_{n+1}=\del_nu$. Similarly, $\del'_{n+1}=u\del'_n$.\square

We will now again consider a cobordism $W$ as in the previous section,
using the same notation.

\be{prop}\label{w-delta}
If $W$ is a cobordism as in Section~\ref{sec:umapcob} then the following hold:
\be{description}
\item[(i)]If $b^+(W)\ge1$ then for all $j\ge0$,
\[\del_jW^*=0,\quad W^*\del'_j=0.\]
\item[(ii)]
If $b_1(W)=b^+(W)=0$ and $k\ge0$ then
\[\be{array}{lll}
\del_jW^*=0,&W^*\del_j'=0 & \text{for $0\le j\le k-1$,}\\
\del_jW^*=\del_{j-k}, & W^*\del_j'=\del_{j-k}' & \text{for $j\ge k$}.
\end{array}\]
\end{description}
\end{prop}

\proof We first introduce some notation that will be used in both parts of the
proof. Choose a base-point $y_j\in Y_j$.
We take $x_j:=(0,y_j)$ as base-point for the band $K_j:=[-1,1]\times Y_j$.
The natural complex line bundle over $\cb^*(K_j)$ will
be denoted $\bl_j$.  For $j=1,2$ choose
$\eps_j>0$ such that there is no monopole $\beta$ over $Y_j$ with
$0<\cd(\beta)-\cd(\theta_{m_j})\le\eps_j$. For elements $\om$ of a moduli
space $M(\hw;\al,\thet2)$ the $\eps_1$--slice on $\rmy_1$
and the
$\eps_2$--slice on $\rpy_2$ (when defined) will be denoted $\tau_1(\om)$ and
$\tau_2(\om)$, respectively. The restriction of $\om$ to the
band $[\tau_j(\om)-1,\tau_j(\om)+1]\times Y_j$ will be denoted $R_j(\om)$.

In both parts of the proof we assume that metric and perturbations on
$\hw$ are chosen so that the conclusions of Proposition~\ref{sec:pert-red} 
hold.

(i) By assumption there is no reducible monopole in any moduli space
over $\hw$.
If $\eps_2$ is sufficiently small then $\tau_2$ is defined on
any moduli space $ M(\hw;\al,\thet2)$ with $\al\in\fl(Y_1,m_1)$.
Let $s$ be a generic section of $\oplus^j\bl_2$. Set
\[M'_\al:=\{\om\in M(\hw;\al,\thet2)\st s(R_2(\om))=0\}.\]
If $\ind(\al)=2m_2-2j-1+d$ then $M'_\al$ consists of a finite
number of oriented points; for irreducible $\al$
let $\bar\del_{W,j}(\al)$ denote the number of these
points counted with sign. If $\ind(\al)$ is one less, then 
since the number of ends of
the $1$--manifold $M'_\al$, counted with sign,
must be zero, we have
\[(\bar\del_jW^\#-\bar\del_{W,j}d)\al=0,\]
hence $\del_jW^*=0$. Similarly one finds that $W^*\del_j'=0$.

(ii) We will focus on the statement about $\del_j$. The case of
$\del_j'$ is similar. We will write $\fl_j=\fl(Y_j,m_j)$ and
$\cf^*_j=\cf^*(Y_j,m_j)$.

Let $\al\in\fl_1^*$ have index $2m_2-2+d=2(m_1+k)-2$.
Counting with sign the ends of the $1$--manifold
$M(\hw;\al,\thet2)$ we find that
\[(\del W^\#-\del_Wd)\al=
\be{cases}
\del\al & \text{if $k=0$,}\\
0 & \text{if $k\ge1$}.
\end{cases}\]
If $k=0$ then $\del_0W^*=\del_0$ in cohomology, and we are done.
From now on assume $k\ge1$, in which case $\del_0W^*=0$.
Propositions~\ref{wu-comm} and \ref{prop:delnun} yield
\[\del_jW^*=\del_0u^jW^*=\del_0W^*u^j=0,\quad 0\le j\le k-1.\]
For $j\ge k$, Proposition~\ref{wu-comm} yields
\[\del_jW^*=\del_ku^{j-k}W^*=\del_kW^*u^{j-k},\quad j\ge k.\]
To complete the proof of the proposition we need only prove the 
following lemma:

\be{lemma}\label{del-k}
$\del_kW^*=\del_0$.
\end{lemma}

\proof Let $\al\in\fl^*_1$ have index $\ge 2m_1-2$. Then the
moduli space $M_\al:=M(\hw;\al,\thet2)$ has dimension $\le 2k+1$.
We are going to cut down $M_\al$, roughly speaking, by the $k$'th power
of the Chern class of a complex line bundle $\bl_\al$ 
to obtain a manifold $\Si_\al$ of dimension $\le1$.
This line bundle may be thought of
as the natural line bundle associated to a base-point at $\infty$ over
the end $\rpy_2$. (Compare the definition of $\del_n$ above).
The proof will involve studying the ends of $\Si_\al$.
In the case when $\al$ has index $2m_1-2$ special care must be
taken with the part of $M_\al$ consisting of
monopoles whose spinors are small over $W$. By the gluing theorem 
in \cite{Fr13},
such monopoles are obtained by gluing elements of
$M(\al,\thet1):=M(\ry_1;\al,\thet1)$
with elements of $M_\theta:=M(\hw;\thet1,\thet2)$ that are
close to the reducible point $\omred$.
We will use $c_1(\bl_\al)^{k-1}$ to cut down the link of $\omred$, a copy of
$\cp{k-1}$, to a point, while the last factor $c_1(\bl_\al)$ will be played off
against the gluing parameter.

Note that $M(\al,\thet1)$ is empty if $\ind(\al)>2m_1-2$ or
$\cd(\al)\le\cd(\thet1)$. In that
case 
\[\inf_{[A,\phi]\in M_\al}\int_W|\phi|>0\]
and the discussion below simplifies significantly.

First some notation:
Let $x_0\in W$ be a base-point for $W$ and $\hw$.
The natural complex line bundle over $M_\al$
be denoted $\bl_\al$. As in \cite{Fr13} we
denote by $R_K(S)$ or $S|_K$ the restriction  of a configuration $S$
(or a gauge equivalence class of such)
to a codimension~$0$ submanifold $K\subset\hw$.

{\bf (I)} This part of the proof
is concerned with the framed moduli space
$M_{\theta,x_0}:=M_{x_0}(\hw;\thet1,\thet2)$ with base-point $x_0$. To prepare
for the application of the gluing theorem we will construct the following:
\be{itemize}
\item a $\U1$--invariant precompact
open neighbourhood $G_0$ of $\omred$ in $M_{\theta,x_0}$,
\item a $\U1$--invariant open subset $V_0\subset\cb_{x_0}(W)$ containing
$R_W(\oline G_0)$,
\item a $\U1$--equivariant smooth map $q_0:V_0\to M_{\theta,x_0}$ such
that $q_0(\om|_W)=\om$ for all $\om\in\oline G_0$.
\end{itemize}

Let $(\ared,0)$ be a
representative for $\omred$
and set $\bred:=\ared|_W$. By \cite[Section~2.5]{Fr13} there is a closed
$\U1$--invariant subspace
\[\Xi\subset L^p_1(W;i\La^1_W\oplus S^+_W)\]
such that 
\[\Xi':=(\bred,0)+\Xi\to\cb_{x_0}(W)\]
is a local diffeomorphism at $(\bred,0)$.
Since the tangent space of $M_{\theta,x_0}$
at $\omred$ can be identified with $\ker D_{\ared}$, it follows by unique
continuation for harmonic spinors that the restriction map
\[R_W:M_{\theta,x_0}\to\cb_{x_0}(W)\]
is immersive at $\omred$. Let $G^+_0\subset M_{\theta,x_0}$ be an open subset
such that
\be{itemize}
\item $G^+_0$ is the image of a $\U1$--equivariant
smooth embedding
\[\iota:\co^k\to M_{\theta,x_0},\]
\item $R_W$ maps $G^+_0$ diffeomorphically onto
a submanifold $H\subset\cb_{x_0}(W)$,
\item there is a $\U1$--invariant submanifold
$\ti H\subset\Xi'$ which maps diffeomorphically onto $H$.
\end{itemize}
Let 
$\Xi_0\subset\Xi$ be a closed linear complement of the tangent
space $T_{(\bred,0)}\ti H$, so that $\Xi=T_{(\bred,0)}\ti H\oplus\Xi_0$. Then
\[\mu:\ti H\times\Xi_0\to\cb_{x_0}(W),\quad(S,s)\mapsto[S+s]\]
is a $\U1$--equivariant local diffeomorphism at $((\bred,0),0)$,
by the inverse function theorem.
If $\pi_1:\ti H\times\Xi_0\to\ti H$ then $\pi_1\circ\mu\inv$ is a
$\U1$--equivariant smooth map $V_0\to\ti H$ for some $\U1$--invariant
open neighbourhood $V_0$ of $\omred|_W$. Let $\chv0$ denote the image
of $V_0$ in $\cb(W)$ and $\chv0^*$ the irreducible part of $\chv0$.
Let $f:\ti H\to G^+_0$ be the obvious map,
such that $[S]=f(S)|_W$ for all $S\in\ti H$. Set
\[q_0:=f\circ\pi_1\circ\mu\inv:V_0\to G^+_0\subset M_{\theta,x_0}\]
and
\[z:=\iota\inv:G^+_0\to\co^k,\]
and let $z_j$ be the $j$'th coordinate of $z$. If $U$ is any open subset of
an unframed moduli space over $\hw$ which is mapped into $\chv0^*$ by $R_W$
then the $\U1$--equivariant map
\be{equation}\label{eqn:sjdef}
z_j\circ q_0:V_0\to\co
\end{equation}
defines a section $\si_j$ of the appropriate complex line bundle over $U$. 

Let $G_0\subset M_{\theta,x_0}$ be the image under $\iota$
of the $\rho$--ball about the origin
in $\co^k$, where $\rho>0$ is small. In particular, the closure
$\oline G_0$ of $G_0$ should be contained in the image of $q_0$,
and $R_W(\oline G_0)\subset V_0$.

{\bf (II)} This part is concerned with the moduli space $M(\al,\thet1)$ and
is only relevant in the case when $\ind(\al)=2m_1-2$ and
$\cd(\al)>\cd(\thet1)$. Then
$M(\al,\thet1)$ has dimension~$1$. 
Again we must construct certain data $G_1,V_1,q_1$ needed for the
application of the gluing theorem.

For any configuration $S$ over $K_1$ set
\[J(S):=\cd(S_0),\]
where $S_0=S|_{\{0\}\times Y_1}$.
Because of the lack of unique continuation for our perturbed equations
we introduce a larger band $K'_1:=[-\nu,\nu]\times Y_1$, where $\nu\ge1$ is
chosen so large that elements of $M(\al,\thet1)$ with $J=\cd(\thet1)+\eps_1$ are
distinguished by their restrictions to $K'_1$.
(The constant $\eps_1$ was chosen at the beginning of the proof of
Proposition~\ref{w-delta}.)
For small $\eps_1'>0$ the set
\[G_1:=\{\om\in M_{x_1}(\al,\thet1)\st|J(\om)-\cd(\thet1)-\eps_1|<\eps_1'\}\]
is precompact in $M_{x_1}(\al,\thet1)$ and one can find
a $\U1$--invariant open neighbourhood $V'_1$ of
$R_{K'_1}(\oline G_1)\subset\cb^*_{x_1}(K'_1)$ and a $\U1$--invariant
open subset $V_1\subset\cb^*_{x_1}(K_1)$ such that
\be{description}
\item[(i)]$R_{K_1}(V'_1)\subset V_1$,
\item[(ii)]if $\chv1\subset\cb(K_1)$ denotes the image of $V_1$,
then $V_1\to\chv1$ is a trivial $\U1$--bundle,
\item[(iii)]for every $\om\in\chv1$ there exists a (necessarily unique)
$\ti\om\in M(\al,\thet1)$ with $J(\om)=J(\ti\om|_{K_1})$,
\item[(iv)]$\chv1'$ is a disjoint union of open sets, one for each component
of $\check G_1$ (the image of $G_1$ in $M(\al,\thet1)$),
such that $R_{K'_1}$ maps each of these components into the
corresponding open subset of $\chv1'$.
\end{description}
The $\ti\om$ in (iii) will be denoted $\check q_1(\om)$.
Choose a section $e$ of the bundle $V_1\to\chv1$, and a section
$e'$ of the bundle $M_{x_1}(\al,\thet1)\to M(\al,\thet1)$ which agrees with
the pull-back
section $(R_{K_1})^*e$ on $\check G_1$.
Let $q_1:V_1\to M_{x_1}(\al,\thet1)$ be the
$\U1$--equivariant map satisfying $q_1\circ e=e'\circ\check q_1$.

{\bf (III)} We will now construct the sections of
$\bl_\al$ that will be used to define $\Si_\al$.
For any configuration $(A,\Phi)$ over $\hw$ set
\be{equation}\label{eqn:elldef}
\ell(A,\Phi)=\int_W|\Phi|.
\end{equation}
If $\eps_0>0$ is sufficiently small then the region 
$U\subset M_\al$ defined by the inequality $\ell<\eps_0$
(which is empty if $M(\al,\thet1)$ is empty) can be described by the
gluing theorem \cite[Theorem 11.1.1]{Fr13} with $G=G_1\times G_0$,
$K=K_1\coprod W$, $V=V_1\times V_0$, and $q=q_1\times q_0$.
(Thus, roughly speaking,
all elements of $U$ can be obtained by gluing elements of $G_1$ with elements
of $G_0$.) We will not spell this out precisely, but contend ourselves 
with a list of the conditions on $\eps_0$ most relevant for the exposition 
(including two conditions not directly related to the gluing theorem):

\be{description}
\item[(v)]for any $\om\in U$ the
$\eps_1$--slice $\tau_1(\om)$ on $\rmy_1$
is defined, $\tau_1(\om)+\nu<0$, and the restriction of $\om$ to the band
$[\tau_1(\om)-\nu,\tau_1(\om)+\nu]\times Y_1$ lies in $\chv1'$,
\item[(vi)]$\om|_W\in\chv0$ for all $\om\in U$,
\item[(vii)]all $\om\in M_\theta$ with $\ell(\om)<2\eps_0$ lie in the image
$\check G_0$ of $G_0$ in $M_\theta$,
\item[(viii)]$\ell>\eps_0$ on any moduli space $M(\hw;\beta_1,\beta_2)$ in
which $\inf(\ell)>0$.
\end{description}

Note that the chain-limit (see \cite[Section~7.1]{Fr13}) of a sequence in $U$
must either lie in $M_\al$ or have the form $(\om_1,\om_0)$, where
$\om_1\in M(\al,\thet1)$ and $\om_0\in\check G_0$. 
The map $U\to\chv1'$ indicated in (v) together with the decomposition
of $\chv1'$ in (iv) gives rise to a map $\psi:U\to\chm(\al,\thet1)$.

Let $\ga:\R\to\hw$ be a smooth path such that $\ga(0)=x_0$ and
\[\ga(t)=
\be{cases}
(t,y_1) & t\le-1,\\
(t,y_2) & t\ge1.
\end{cases}\]
Let $\cb_x:=\cb(\hw;\al,\thet2)$ be the framed orbit space of configurations
with base-point $x$ (see \cite[Section~3.4]{Fr13}).
For $a<b$ we define the holonomy map
\[h_a^b:\cb_{\ga(a)}\to\cb_{\ga(b)}\]
(depending on the choice of a reference spin connection $A\rr$ over $\hw$)
as follows: If $S=(A,\Phi)$ is a representative of $\om\in\cb_{\ga(a)}$ and
$u:\hw\to\U1$ an element of the group of gauge transformations acting
on the configuration space $\cc(\hw;\al,\thet2)$ such that
\[u(\ga(b))=\exp\left(\int_a^b\ga^*(A-A\rr)\right),\]
then
$u(S)$ is a representative of $h_a^b(\om)$. (Cf.\ \cite[Eqn.~11.1]{Fr13}.)
One can easily verify that this holonomy map behaves functorially with
respect to composition of paths.

That the map $R_1:U\to\cb^*(K_1)$ (defined in the very beginning of the
proof of Proposition~\ref{w-delta})
is smooth follows from \cite[Proposition~3.4.2 and Lemma~3.3.1]{Fr13}
By means of holonomy along $\ga$ the pull-back
$R_1^*(e)$ defines a smooth section $\si^-$ (of constant norm $1$) of
$\bl_\al|_U$. Because of the condition (vi) above the map \Ref{eqn:sjdef}
provides sections $\si_1,\dots,\si_k$ of $\bl_\al|_U$.
On the part of $U$ where
$\si_1\neq0$ we define a map $\eta$ into $\U1$ by the condition
\[\si^-=\eta\cdot\si_1/|\si_1|.\]

Choose a smooth function $c:\R\to(1,\infty)$ such that $c(t)=t$ for
$t\ge4$ and $c(t)=2$ for $t\le3$. For $\om\in M_\al$ set
\be{equation}\label{tau-cut-off}
\bar\tau_2(\om)=
\be{cases}
c(\tau_2(\om)) & \text{if $\tau_2(\om)$ is defined},\\
2 & \text{otherwise}.
\end{cases}
\end{equation}
For any $\om\in M_\al$ let $\bar R_2(\om)\in\cb^*(K_2)$
denote the restriction of
$\om$ to the band $[\bar\tau_2(\om)-1,\bar\tau_2(\om)+1]\times Y_2$.
Choose generic smooth sections $s'_1,\dots,s'_k$ of $\bl_2$.
By means of holonomy along $\ga$
the pull-back $(\bar R_2)^*s'_j$ defines a section $\si'_j$ of $\bl_\al$. 
In a similar way one obtains a section (also denoted $\si'_j$) of the natural
complex line bundle over $M^*_\theta$.

Choose a smooth map $\xi:\R\times\U1\to\co$ such that
\[\xi(t,z)=
\be{cases}
1,&t\ge2,\\
z,&t\le1,
\end{cases}\]
and such that $0$ is a regular value of $\xi$. Then $\#\xi\inv(0)$
is independent of $\xi$, and a simple computation shows that
$\#\xi\inv(0)=-1$. 

Define a smooth section $\ti\si_1$ of $\bl_\al|_U$ as follows:
\[\ti\si_1:=\be{cases}
\si^- & \text{where $\si_1=0$,} \\
\xi(-\tau_1|\si_1|,\eta)\cdot\si_1/|\si_1| & \text{where $\si_1\neq0$.}
\end{cases}\]
Note that if $\om\in U$ and $\si_1(\om)=0$ then $\ti\si_1=\si^-$
in a neighbourhood of $\om$. To see this, observe that $\tau_1$ is bounded
in any compact neighbourhood of $\om$, hence
$-\tau_1|\si_1|\le1$ in a small neighbourhood. 
On the other hand, if $|\si_1(\om)|\ge2$ then $\ti\si_1=\si_1/|\si_1|$ at $\om$,
since $\tau_1<-1$ by definition.

Choose a smooth function $\ti\ka:\R\to\R$ such that
$\ti\ka(t)=1$ for $t\ge\eps_0/2$,
$\ti\ka(t)>0$ for $t>\eps_0/4$, and $\ti\ka(t)=0$ for $t\le\eps_0/4$. Set
\[\ka:=\ti\ka\circ\ell,\quad\uka:=\ti\ka\circ(\ell/2),\]
where $\ell$ is as in \Ref{eqn:elldef},
and define smooth sections
$\hat\si_1,\dots,\hat\si_k$ of $\bl_\al$ by
\be{align*}
\hat\si_1&=\ka\si'_1+(1-\ka)\ti\si_1,\\
\hat\si_j&=\uka\si'_j+(1-\uka)\si_j,\quad 2\le j\le k.
\end{align*}
Note that the interpolation regions are different in the two cases:
For $\hat\si_1$ the region is given by $\eps_0/4<\ell<\eps_0/2$, whereas for
$\hat\si_j$, $j>1$ it is $\eps_0/2<\ell<\eps_0$. This will be essential
in the compactness argument in part~(V), Case~1 below.

{\bf (IV)} We will now derive a gluing result which will be used in (V).
For any $\tau>0$ let $\hut$ be the part of $U$ where $-\tau_1>\tau$,
$0<|\si_1|<\rho$, and $\si_j=0$ for $j=2,\dots,k$. Applying the
gluing theorem \cite[Theorem~14.3.1]{Fr13} with $G=G_1\times G_0$,
$K=K_1\coprod W$, $V=V_1\times V_0$, and $q=q_1\times q_0$ we find that
if $\tau\gg0$ then $\hut$ is a smooth submanifold of $M_\al$ and the map
\[(\eta,-\tau_1,|\si_1|,\psi):\hut\to\U1\times(\tau,\infty)\times
(0,\rho)\times\chm(\al,\thet1)\]
is an orientation reversing diffeomorphism. 
(Argue by contradiction to prove the map is injective for $\tau\gg0$.
Another application of the gluing theorem shows that the map is also a 
surjective submersion for large $\tau$. To obtain submersiveness it is useful
to vary the parameter $T$ in the gluing theorem as well.)

{\bf (V)} Let $\Si_\al\subset M_\al$
be the locus where $\hat\si_j=0$ for $j=1,\dots,k$. If $\eps_0>0$ is 
sufficiently small then $\Si_\al$ is transversely cut out from $M_\al$.
If $\al$ has index $2m_1-1$ then
$\Si_\al$ consists of a finite number of points,
and we define
\[\del_{W,k}:\cf_1^{2m_1-1}\to\z,\quad\beta\mapsto\#\Si_\beta.\]

Now let $\al$ have index~$2m_1-2$, in which case $\Si_\al$ has dimension~$1$.
We will determine the ends of $\Si_\al$.
Let $\{\om_n\}$ be a sequence in $\Si_\al$.
After passing to a subsequence we may assume that $\{\om_n\}$
chain-converges to some triple consisting of a broken gradient line over
$\R\times Y_1$, a monopole $\om$ over $\hw$, and a broken gradient
line over $\R\times Y_2$. Suppose $\om\not\in M_\al$ (or equivalently,
$\om\not\in\Si_\al$).

{\bf Case 1:} $\om\in M_\theta$. We will show shortly that
$\om$ must be the reducible point $\omred$. Granted this we deduce from
(IV) that the corresponding 
number of ends of $\Si_\al$ is $(\#\xi\inv(0))\del\al=-\del\al$.

Suppose now to the contrary that $\om\neq\omred$. If
$\ell(\om)>\eps_0/2$ then $\om$ would satisfy the equations
\be{gather*}
\si'_1=0,\\
\uka\si'_j+(1-\uka)\si_j=0,\quad j=2,\dots,k.
\end{gather*}
However, the zero-set of these equations in $M^*_\theta$ is empty
for a generic choice of sections $s'_j$, because $\dim M^*_\theta=2k-1$.
Hence $\ell(\om)\le\eps_0/2$, and consequently $\om\in\check G_0$ by
condition~(vii) above.
Certainly, $\om$ must satisfy the equations
\be{equation}\label{eqn:siom1}
\si_j=0,\quad j=2,\dots,k.
\end{equation}
Therefore $\si_1(\om)\neq0$. Since $\tau_1(\om_n)\to-\infty$ as $n\to\infty$,
we conclude that 
\[\ti\si_1(\om_n)=\frac{\si_1(\om_n)}{|\si_1(\om_n)|}\quad\text{for $n\gg0$,}\]
so $\om$ satisfies the equation
\be{equation}\label{eqn:siom2}
\ka\si'_1+\frac{1-\ka}{|\si_1|}\,\si_1=0.
\end{equation}
However, for dimensional reasons the common zero-set of the equations
\Ref{eqn:siom1} and \Ref{eqn:siom2} in $\check G_0\setmin\{\omred\}$
is empty for a generic choice of section $s'_1$. This is a contradiction,
so we conclude that $\om=\omred$ as claimed.

{\bf Case 2:} $\om\not\in M_\theta$. Then $\ell(\om)>\eps_0$ by
condition~(viii) above, so for large $n$ one has
$\si'_j(\om_n)=0$, $j=1,\dots,k$.
There are now two possibilities: If factorization
over $\R_+\times Y_2$ does not occur then 
$\bar\tau_2(\om_n)\to\bar\tau_2(\om)$ and $\om$ will lie in some $2k$--dimensional
space $M(\hw;\beta_1,\thet2)$ where $\beta_1\in\fl^*_1$.
The corresponding number of ends of $\Si_\al$
is $-\bar\del_{W,k}d\al$. If such factorization does occur then
$\tau_2(\om_n)\to\infty$ and $\om$ must lie in some $0$--dimensional space
$M(\hw;\al,\beta_2)$ where $\beta_2\in\fl^*_2$. 
The corresponding number of ends of $\Si_\al$ is $\bar\del_kW^\#\al$.

Since the total number of ends of $\Si_\al$ must be zero, we conclude, after
passing to Floer cohomology, that $\del_0=\del_kW^*$.
This completes the proof of the lemma, and also of
Proposition~\ref{w-delta}.\square

\section{Irreducible Floer groups for different chambers}
\label{sec:diff-chambers}

Let $Y$ be a $\spc$ rational homology $3$--sphere, $m\in\mm(Y)$, and
$G$ an Abelian group. In this section we will determine the kernel and image
of the canonical homomorphism
\be{equation}\label{eqn:mtomell}
\hf^*(Y,G;m)\to\hf^*(Y,G;m+\ell)
\end{equation}
for any natural number $\ell$. We first consider the case $\ell=1$ and
then do the general case by induction on $\ell$. For $\ell=1$
Proposition~\ref{w-delta}~(ii) says that the map \Ref{eqn:mtomell} induces
a  homomorphism
\[\bar J:\coker(\del'_0)\to\ker(\del_0).\]

\be{prop}\label{bar-j-iso}
$\bar J$ is an isomorphism.
\end{prop}

The proof is divided into three lemmas.
We begin by constructing two graded homomorphisms
\[Q_1,Q_2:\cf^*(Y,m+1)\to\cf^*(Y,m).\]
such that $Q=Q_1+Q_2$ will induce an inverse of $\bar J$.

Fix a point $y_0\in Y$. Let $D\subset(-1,1)\times Y$ be a compact
$4$--ball disjoint from $(-1,1)\times\{y_0\}$, and set
$W=[-1,1]\times Y\setminus\text{int}(D)$. We regard $\hw$ as a manifold
with two negative ends $\rmy$ and $\rms$ (in that order) and one positive end
$\rpy$.
The notation $M(\hw;\al_-,\al_+,\beta)$ will refer to a moduli space
over $\hw$ with limits $\al_\pm$ over $\R_\pm\times Y$ and $\beta$ over
$\rms$. Note that
\[\dim M(\hw;\al_-,\al_+,\beta)=
-\ind(\al_-)+(\ind(\al_+)-n_{\al_+})-\ind(\beta).\]

If instead one regards $\hw$ as a manifold with one negative end $\rmy$
and two positive ends $\R_+\times S^3$ and $\rpy$, then one
obtains in general a different orientation of $M=M(\hw;\al_-,\al_+,\beta)$. If 
$\al_+$ and $\beta$ are both irreducible then the two orientations
differ by the sign $(-1)^{\ind_2(\beta)(d+\ind_2(\al_+))}$, where
$d\equiv\dim M\mod2$.

For the construction of $Q$ we use a GM-pair for $\hw$ of index
$m+1$, $m$, $-1$ on the ends $\rmy$, $\rpy$, $\R_-\times S^3$, resp.
An important point in the proof of the proposition will be that
\[\dim M(\hw;\theta_{m+1},\theta_m,\theta_{-1})=-1.\] 
Generically, this moduli space consists of the reducible monopole $\omred$
only, 
and this is a regular point.

Fix a small $\eps>0$ and for any element $\om$ of a moduli space
$M(\hw;\al,\beta,\theta_{-1})$ where $\al,\beta$ are not
both reducible
let $\tau(\om)$ be the $(-\eps)$--slice
over the end $\rms$, when this is defined. Let $c$ be as in \Ref{tau-cut-off}
and define
\[\bar\tau(\om)=
\be{cases}
-c(-\tau(\om)) & \text{if $\tau(\om)$ is defined},\\
-2 & \text{otherwise}.
\end{cases}\]
Set $K=[-1,1]\times S^3$ and let $R(\om)\in\cb^*(K)$ denote the
restriction of $\om$ to the band $[\bar\tau(\om)-1,\bar\tau(\om)+1]\times S^3$.
Choose a generic smooth section $s$ of the
natural complex line bundle over $\cb^*(K)$ and set
\[M'(\al,\beta)
:=\{\om\in M(\hw;\al,\beta,\theta_{-1})\st s(R(\om))=0\}.\]
If $\al,\beta$ are both irreducible and of the same index then
$M'(\al,\beta)$ is a finite set of points, so we can define $Q_1$ on
generators $\al$ by
\[Q_1\al=\sum_\beta(\#M'(\al,\beta))\beta.\]
According to  Proposition~\ref{s3vanish} we have $\hf^1(S^3,-1)=0$, so
there is a cochain $z\in\cf^0(S^3,-1)$ such that $\bar\del'_1+dz=0$. 
Define $Q_2$ by
\[Q_2\al=\sum_\beta(\#M(\hw;\al,\beta,z))\beta,\]
where we regard $M(\hw;\al,\beta,z)$ formally as a linear function in $z$.
(To simplify language, we will often speak of such formal linear combinations
of moduli spaces as a single moduli space.)
Set $Q=Q_1+Q_2$.

\be{lemma}\label{dRRd}
$dQ-Qd+\del'\del=0$.
\end{lemma}
Here $\del:\cf^{2m}(Y,m+1)\to\z$ and $\del':\z\to\cf^{2m+1}(Y,m)$ are the
maps defined in Section~\ref{sec:irrflg}.

\proof Let $\al\in\fl^*(Y,m+1)$ and $\beta\in\fl^*(Y,m)$ satisfy
$\ind(\beta)=\ind(\al)+1$. Choose $\eps_1>0$ such that there is no monopole
$\ga\in\fl^*(Y,m+1)$ with $0<\cd(\ga)-\cd(\theta_{m+1})\le\eps_1$.
Choose $\eps_2>0$ such that there is no monopole
$\ga\in\fl^*(Y,m)$ with $0<\cd(\theta_m)-\cd(\ga)\le\eps_2$. For any
configuration $(A,\Phi)$ over $\hw$ set
\be{equation}\label{eqn:ellaphi}
\ell(A,\Phi):=\int_W|\Phi|.
\end{equation}
Note that if $\{\om_n\}$ is any chain-convergent
sequence in $M(\hw;\al,\beta)$ with $\ell(\om_n)\to0$ then the chain-limit
has the form $(\om',\omred,\om'')$, where $\om'\in \check M(\al,\theta_{m+1})$,
$\om''\in\check M(\theta_m,\beta)$.
Fix a small $\eps_0>0$. If $\om\in M(\hw;\al,\beta)$ and $\ell(\om)<\eps_0$
then the $\eps_1$--slice of $\om$ on $\rmy$ and the $\eps_2$--slice
of $\om$ on $\rpy$ will be well-defined; we denote these by $\tau_1(\om)$
and $\tau_2(\om)$, resp. Let $R_j(\om)$ denote the restriction of $\om$
to the band $[\tau_j(\om)-1,\tau_j(\om)+1]$.

Let $\bl\to M(\hw;\al,\beta)$ be the 
natural complex line bundle associated to the base-point $(0,y_0)$.
Let $q:(-\infty,1]\to\hw$
be a path with $q(1)=(0,y_0)\in W$ and $q(t)=(t,y_{-1})$ for $t\le0$, where 
$y_{-1}\in S^3$.
We use holonomy along the paths $\R\times\{y_0\}$ and $q$ to 
identify $\bl$ with the natural complex line bundles associated to other
base-points along these paths. Much as in the proof of Lemma~\ref{del-k} we 
construct a section $\hat\si$ of $\bl$ such that
\be{itemize}
\item $\hat\si$ agrees with the pull-back
section $R^*s$ in the region where $\ell\ge\eps_0/2$,
\item For $j=1,2$, if $\ell(\om)\le\eps_0/4$ and
$(-1)^j(\tau_1(\om)+\tau_2(\om))\ge1$ then
$\hat\si(\om)$ is non-zero and depends only on the restriction of $\om$
to $[\tau_j(\om)-1,\tau_j(\om)+1]\times Y$.
\end{itemize}
In the intermediate regions we interpolate as in the proof of
Lemma~\ref{del-k}.
Let $\Si\subset M(\hw;\al,\beta,\theta_{-1})$ be the zero-set of $\hat\si$.

Now let $\{\om_n\}$ be a chain-convergent sequence
in $\Si$ whose limit does not lie in $\Si$.
If $\bar\tau(\om_n)\to\infty$ then factorization takes place over
$\R_-\times S^3$ through some $\ga\in\fl^*(S^3,-1)$ of index~$1$;
the corresponding number of ends of $\Si$ (counted with sign)
is $\la\chi,\bar\del'_1\ra$, where
\[\chi=\sum_\ga(\#M(\hw;\al,\beta,\ga))\ga\in\cf_1(S^3,-1).\]
If $\bar\tau(\om_n)$ stays bounded, then there are two cases depending on
whether the chain-limit involves the reducible monopole over $\hw$ or not.
In the reducible case one finds as in the proof of Proposition~\ref{dvvd}
that the corresponding number of ends of $\Si$ is $\la\del'\del\al,\beta\ra$.
In the irreducible case factorization takes place over one of the tubular 
ends $\R_\pm\times Y$; here the corresponding number of ends of $\Si$ is
$\la(dQ_1-Q_1d)\al,\beta\ra$. 
Since the total number of ends of $\Si$ must be zero, we conclude that
\[\la(dQ_1-Q_1d+\del'\del)\al,\beta\ra+\la\chi,\bar\del'_1\ra=0.\]
Similarly, by considering the ends of $M(\hw;\al,\beta,z)$ we find that
\[\la(dQ_2-Q_2d)\al,\beta\ra+\la\chi,dz\ra=0.\]
Summing the last two equations and recalling that $\bar\del'_1+dz=0$
we obtain
\[\la(dQ-Qd+\del'\del)\al,\beta\ra.\square\]

Consider an arbitrary $m'\in\mm(Y)$ and
let $B:\cf^*(Y,m')\to\cf^*(Y,m')$ be the homomorphism
defined in the same way as $Q$, but now using a GM-pair for $\hw$ of index
$m'$, $-1$ on the ends $\R_\pm\times Y$, $\R_-\times S^3$, resp.,
and with the same element $z$. Arguing
as in the proof of Lemma~\ref{dRRd} we find that $B$ is a cochain map
(there are no complications from reducibles here). 

\be{lemma}\label{B-id}
$B$ induces the identity map on $\hf^*(Y,m';G)$.
\end{lemma}

As there are two signs in the proof which the author has
not computed (and which affect the sign of $I$ in formula~\Ref{eqn:NI}),
this lemma as well as the next one really hold only up to sign.
However, we will ignore this, as it is irrelevant to the proof of
Proposition~\ref{bar-j-iso}.

\proof  Set $Z:=[-1,0]\times S^3$, and let $\smo$ and $\sze$ be two copies
of the three-sphere. We will think of $\hz$ as having tubular ends
$\R_-\times\smo$ and $\R_+\times\sze$, and the end $\R_-\times S^3$ of $\hw$
will now
be denoted $\R_-\times\sze$. Choose a GM-pair for
$\hw$ which has index $m'$ on each of the ends $\R_\pm\times Y$ and
which agrees with $(1\times\ti g,0)$ on $\R_-\times\sze$, where $\ti g$ is
the round metric on the three-sphere. Choose a GM-pair for $\hz$ which has
index $-1$ on the end $\R_-\times\smo$ and which 
agrees with $(1\times\ti g,0)$ on $\R_+\times\sze$. Let $\xt$ denote the
manifold with tubular ends obtained from $\hw\cup\hz$ by replacing the
ends $\R_\pm\times\sze$ with a neck $[-T,T]\times\sze$. Let $\xt$ be equipped
with the GM-pair inherited from the GM-pairs of $\hw,\hz$
in the natural way.

Now fix
$\al,\beta\in\fl^*(Y,m')$ of the same index and for any $\ga\in\fl(S^3,-1)$
define the parametrized moduli space
\[\cm_\ga:=\bigcup_{T\ge1}M(\xt;\al,\beta,\ga)\times\{T\}.\]
Set $\cm:=\cm_{\theta_{-1}}$. Fix a small $\eps_0>0$ and 
let $U$ be the set of all $(\om,T)\in\cm$ with $\ell(\om)<\eps_0$, where
\[\ell(A,\Phi):=\int_Z|\Phi|.\]
We now copy the construction in the proof of Lemma~\ref{del-k}
(with $k=1$) of a section
$\hat\si_1$ of the natural complex line bundle over 
$\cm$ corresponding to some base-point in $Z$. Here the present
$Z$, $M(\hw;\al,\beta,\theta_0)$, $-T$ correspond to the
$W$, $\check M(\al,\thet1)$, $\tau_1$ resp.\ in the proof of Lemma~\ref{del-k}.
Let $\cm'\subset\cm$ be the zero-set of $\hat\si_1$ and set
\[\cm'':=\cm'+\cm_z,\]
regarded as a formal linear combination of oriented, smooth $1$--manifolds with
boundary.
With the orientation conventions in \cite{Fr13} one has
\[\#\prtl\cm''=-\la B\al,\beta\ra.\]
We will now describe the ends of $\cm''$. Because
$\dim M(\hz;z,\theta_0)=-1$ and there are no non-empty moduli spaces
over $\hz$ of negative expected dimension, it follows that
$M(\xt;\al,\beta,z)$ is empty for large $T$. 
Therefore the ends of $\cm_z$ correspond to
factorization through irreducible critical points
over the tubular ends $\R_\pm\times Y$ or $\R_-\times S^3_{-1}$, which may
occur for finitely may values of $T$.

Determining the ends of $\cm'$ requires a little more consideration. Let 
$\{(\om_n,T_n)\}$ be a sequence in $\cm'$.
By passing to a subsequence we may arrange that the
sequence chain-converges. This means in particular that
$\{T_n\}$ converges to some point $T_\infty\in[1,\infty]$.
Suppose the chain-limit does not lie in $\cm$ (or equivalently, not in $\cm'$).
We consider first the case $T_\infty<\infty$. Observe that
\[\dim M(X^{(T_\infty)};\theta_{m'},\theta_{m'},\theta_{-1})=1.\]
Therefore, for
dimensional reasons the limit over $X^{(T_\infty)}$ cannot be reducible.
Hence factorization through an irreducible critical point has taken place
over $\R_\pm\times Y$ or $\R_-\times S^3_{-1}$. 
We now turn to the case $T_\infty=\infty$. Since
\[\dim M(\hz;\theta_{-1},\theta_0)=1,\]
the chain-limit must have the form
$(\om,\om')$, where $\om\in M(\hw;\al,\beta,\theta_0)$ and
$\om'\in M(\hz;\theta_{-1},\theta_0)$. As in 
the proof of Lemma~\ref{del-k} (Part~(V), Case~1)
one finds that $\om'$ must be the
reducible point and that the corresponding number  $N$ of ends
of $\cm'$ is equal to $\#M(\hw;\al,\beta,\theta_0)$.
Let $W\tu$ be the result of capping off the end $\R_-\times S^3_0\subset\hw$
with a disk, allowing for a long neck
$[-T,T]\times S^3_0$. An easy version of the argument in
\cite[Part~3]{Fr13} yields
\[\#M(\hw;\al,\beta,\theta_0)=\#M(W\tu;\al,\beta).\]
Deforming the metric on $W\tu$ back to a cylindrical one we conclude that
\be{equation}\label{eqn:NI}
N=\la(I+f_2d\pm df_2)\al,\beta\ra
\end{equation}
for some homomorphism $f_2:\cf^*(Y,m')\to\cf^{*-1}(Y,m')$.

Summing up, recall that the number of boundary points plus the number
of ends of $\cm''$ must be zero, and note that the number of
ends resulting from factorizations over $\R_-\times S^3_{-1}$ cancel out
as in the proof of Lemma~\ref{dRRd}. It follows that
\[\la(-B+I+fd\pm df)\al,\beta\ra=0,\]
where $f$ is the sum of $f_2$ and a homomorphism of degree $-1$ related to
factorizations over $\R_\pm\times Y$.\square

Lemma~\ref{dRRd} shows that $Q$ induces a homomorphism
\[\bar Q:\ker(\del_0)\to\coker(\del'_0).\]

\be{lemma}\label{JR}
$\bar J$ and $\bar Q$ are inverse maps.
\end{lemma}

\proof We show that $\bar J\bar Q=\id$, the proof of $\bar Q\bar J=\id$
being similar.
Set $Z:=[0,1]\times Y$. Choose a GM-pair for $\hw$ which has index
$m+1,m,-1$ over the ends $\rmy$, $\rpy$, $\rms$, resp., and a GM-pair for
$\hz$ which has index $m,m+1$ over the ends $\rmy$, $\rpy$, resp., and such
that over the ends $\rpy\subset\hw$ and $\rmy\subset\hz$ these
are induced from the same ND-pair for $Y$. Let $\xt$ denote the
manifold with tubular ends obtained from $\hw\cup\hz$ by replacing the
ends $\rpy\subset\hw$ and $\rmy\subset\hz$ with a neck $[-T,T]\times Y$.
Let $\xt$ be equipped with the GM-pair inherited from those of $\hw,\hz$.

Now fix $\al,\beta\in\fl^*(Y,m+1)$ of the same index
and define the parametrized moduli spaces $\cm_\ga$ and $\cm$ 
as in the proof of Lemma~\ref{B-id}. Let $\ell$ be as in \Ref{eqn:ellaphi}.
Note that if $\{(\om_n,T_n)\}$ is any chain-convergent sequence in $\cm$
with $\ell(\om_n)\to0$ then $T_n\to\infty$ and the chain-limit is made up of
the reducible monopole over $\hw$ and one element from each of the
sets $\check M(\al,\theta_{m+1})$ and $M(\hz;\theta_m,\beta)$ (which must then 
both be finite). We now copy the construction of $\hat\si$ in the proof of
Lemma~\ref{dRRd}, with $T$ and $M(\hz;\theta_m,\beta)$ playing the roles
of $\tau_2$ and $\check M(\ry;\theta_m,\beta)$, resp. Let $\cm'\subset\cm$
be the zero-set of $\hat\si$ and set
\[\cm'':=\cm'+\cm_z.\]
By considering the boundary points and ends of $\cm''$
one finds that
\be{equation}\label{JRB}
JQ=B+\hat\del'\del+df\pm fd,
\end{equation}
where $J:\cf^*(Y,m)\to\cf^*(Y,m+1)$ is a cochain map which induces the
canonical map in cohomology,
$B$ is of the kind discussed in Lemma~\ref{B-id}, and
\[\hat\del'(1)=\sum_\ga(\#M(\hz;\theta_m,\ga))\ga\in\cf^{2m}(Y,m+1).\]
Combining equation~\Ref{JRB} with Lemma~\ref{B-id} we conclude that
$\bar J\bar Q=\id$.
This completes the proof of Lemma~\ref{JR}, and
also of Proposition~\ref{bar-j-iso}.\square


\be{prop}\label{prop:varchamb}
Let $\ell$ be any natural number. Then the kernel and image of the
natural homomorphism
\[\hf^*(Y,G;m)\to\hf^*(Y,G;m+\ell)\]
are
\[\sum_{j=0}^{\ell-1}\im(\del'_j)
\qquad\text{and}\qquad\bigcap_{j=0}^{\ell-1}\ker(\del_j),\]
respectively.
\end{prop}

\proof For both the kernel and the image, one inclusion follows from
Proposition~\ref{w-delta}~(ii), whereas the other inclusion is easily proved by
induction on $\ell$ using Propositions \ref{w-delta}, \ref{bar-j-iso} and
the functoriality of the homomorphisms
in question with respect to composition of cobordisms.\square

We will now define the maps $D,D'$ in the sequence \Ref{eqn:exact-seq}
and prove that the sequence is exact. 
Fix $m\in\mm(Y)$ and let 
\be{gather*}
v:\uhf^{2m+1}(Y;G)\oset\approx\longrightarrow\hf^{2m+1}(Y,m;G),\\
w:\hf^{2m}(Y,m+1;G)\oset\approx\longrightarrow\ohf^{2m}(Y;G)
\end{gather*}
be the canonical maps. For $z\in\ohf^{2m}(Y;G)$ and $g\in G$ we define
\[D(z):=x^m\otimes\del_0w\inv(z),\qquad D'(x^m\otimes g):=v\inv\del'_0(g).\]
We set $D=0$ on $\ohf^{2m+1}(Y;G)$. (There is little choice,
since $P^{2m+1}(Y)=0$.)

From Proposition~\ref{w-delta} we deduce that $D,D'$ are $\z[u]$--homomorphisms
and that the two right-most squares of Diagram~\Ref{diagr:W} commute.

\be{prop}
The fundamental sequence \Ref{eqn:exact-seq} is exact.
\end{prop}

\proof Exactness at the terms $\uhf$ and $\ohf$
follows from Proposition~\ref{prop:varchamb}, and 
Lemma~\ref{dRRd} yields $D'D=0$. It remains to show that
$\ker(D')\subset\im(D)$. Suppose $D'(x^m\otimes g)=0$. This means there is
a chain $\si=\sum_j\al_j\otimes g_j$ with
$\al_j\in\fl(Y,m)$ and $g_j\in G$ such that
\[\rho\otimes g+d\si=0,\]
where $\rho=\del'(1)\in\cf^{2m+1}(Y,m)$.
The oriented $1$--manifold $M(\ry;\theta_m,\theta_{m+1})$
has (generically) exactly one boundary point (namely the reducible one).
Counting with sign we have that the number of 
boundary points plus the number of ends of that moduli space
must be zero, and similarly the number of ends of $M(\ry;\al_j,\theta_{m+1})$
is zero. Thus
\[\pm1=\la \rho,\phi\ra+\del\psi,\qquad
0=\la d\al_j,\phi\ra+\del\psi_j\]
for a suitable chain $\phi\in\cf_{2m+1}(Y,m)$ and cochains
$\psi,\psi_j\in\cf^{2m}(Y,m+1)$.
Tensoring the second equation with $g_j$, summing over $j$, and adding
the first equation tensored with $g$ we obtain $\pm g=\del\tau$, where
\[\tau=\psi\otimes g+\sum_j\psi_j\otimes g_j.\]
A similar argument shows that $d\tau=0$, hence $x^m\otimes g\in\im(D)$.
\square


\section{Reduced Floer groups}

In this section we will show how the reduced Floer cohomology $\rhf^*(Y;G)$
can be computed in an arbitrary chamber $m$ by factoring out interaction in
$\hf^*(Y,m;G)$ with the reducible critical point. We will then prove 
Theorem~\ref{thm:swltr}.

\be{defn}
We define the {\em reduced} Floer cohomology group $\rhf^*(Y,m;G)$ by
\[\rhf^q(Y,m;G)=Z^q/B^q,\]
where
\[Z^*=\bigcap_{j\ge0}\ker(\del_j),\qquad B^*=\sum_{j\ge0}\im(\del'_j).\]
\end{defn}

We define
\[u:\rhf^*(Y,m;G)\to\rhf^{*+2}(Y,m;G)\]
as follows: Choose a map
$v$ as in Proposition~\ref{dvvd}. Let $Z'^*$
(resp.\ $B'^*$) be the set of cocycles in $\cf^*(Y,m)\otimes G$ representing 
an element of $Z^*$ (resp.\ $B^*$). By Proposition~\ref{dvvd} the
endomorphism $v\otimes1$ maps each
of $Z^*$ and $B^*$ into itself and therefore
induces the desired endomorphism of $Z'^*/B'^*=\rhf^*(Y,m;G)$.

Let $W$ be a cobordism as in Section~\ref{sec:umapcob}. If the assumptions
of either part~(i) or part~(ii) of Proposition~\ref{w-delta} are satisfied then
by Propositions~\ref{wu-comm}, \ref{w-delta} the map $W^*$
gives rise to a homomorphism
\[\rhf^*(Y_1,m_1;G)\to\rhf^{*-d}(Y_2,m_2;G)\]
which commutes with the $u$--maps. 

\be{prop}\label{red-iso}
(i) If $m_1\le m_2$ then the $\spc$ cobordism $[0,1]\times Y$
induces an isomorphism
\[\rhf^*(Y,m_1;G)\oset\approx\to\rhf^*(Y,m_2;G).\]
(ii) $\rhf^*(Y,m;G)$ is canonically isomorphic to $\rhf^*(Y;G)$
as a $\z[u]$ module.
\end{prop}

\proof Statement (i) is an immediate consequence of
Propositions~\ref{w-delta}~(ii) and \ref{prop:varchamb}. 
To prove (ii), recall that if
$2m_1\le q\le2m_2-1$ then by Proposition~\ref{prop:irred-iso}
we have a sequence of homomorphisms
\[\uhf^q\oset\approx\to\hf^q_{m_1}\to\hf^q_{m_2}\oset\approx\to\ohf^q,\]
where we have used the same abbreviations as in \Ref{eqn:equmap}.
By Proposition~\ref{prop:varchamb} the image of the middle homomorphism
is precisely $\rhf^q(Y,m_2;G)$. On the other hand, that image can clearly
be identified with $\rhf^*(Y;G)$.\square

{\em Proof of Theorem~\ref{thm:swltr}:} Choose $m\in\mm(Y)$ and a sequence
$\{s_j\}$ of generic sections of the natural complex line bundle over
$\cb^*(W)$. For any $\al,\beta\in\fl^*(Y,m)$ and non-negative integer $k$ let
\[M_{\al\beta k}:=\{\om\in M(\hw;\al,\beta)\st s_j(\om|_W)=0,\quad
j=1,\dots,k\}.\]
Define a cochain map
\[P_k:\cf^*(Y,m)\to\cf^{*-2n+2k}(Y,m)\]
on generators by
\[P_k(\al)=\sum_\beta(\#M_{\al\beta k})\beta.\]
Let $\psi_k$ denote the induced endomorphism of $\hf^*(Y,m)$ (with integral
coefficients). In analogy with Proposition~\ref{w-delta}~(i) one has
\[\del_j\psi_k=0,\quad\psi_k\del'_j=0\]
for all $j\ge0$, hence $\psi_k$ gives rise to an endomorphism 
$\hat\psi_k$ of $\rhf^*(Y,m)$. Choose a map $v$ as in Proposition~\ref{dvvd}. Then
\[P_{k+1}-vP_k+\del'\del_{(k)}=d\phi_k\pm\phi_kd,\]
where $\del_{(k)}$ and $\phi_k$ are maps similar to $\del_W$ and $\phi$
in Proposition~\ref{wu-comm}. By induction on $k$ we conclude that
$\hat\psi_k=u^k\hat\psi_0$ for all $k$. By definition, $\hat\psi_0=\rpsi(W)$.
Now
\[\SW(X)=L(P_n)=L(\psi_n)=L(\psi_n|_{Z^*})=L(\hat\psi_n)=L(u^n\rpsi(W)).\]
The first equality is basic gluing theory, the signs being given by
Proposition~\ref{prop:glue-ortn2}. The second equality is the general fact that
the Lefschetz number of a degree~$0$ chain endomorphism
of a finitely generated chain complex agrees with the Lefschetz number of
the induced endomorphism on homology. The third and fourth equalities
hold because $\psi_n$ takes values in $Z^*$ and vanishes on $B^*$.\square


\section{The $h$--invariant}\label{sec:h-inv}

In this section we will work with Floer cohomology with
coefficients in a fixed field $\f$, which
will usually be suppressed from notation. $Y$ will denote a $\spc$ rational
homology $3$--sphere. All $4$--manifolds will be smooth.

\be{prop}
Let $m\in\mm(Y)$. Then for the group $\hf^*(Y,m)$,
either $\del_0=0$ or $\del'_0=0$.
\end{prop}

\proof This follows immediately from Proposition~\ref{dvvd}.\square

\be{defn}For any $m\in\mm(Y)$ set
\[\zeta=\zeta(Y,m):=\chi(\hf^*(Y,m))-\chi(\rhf^*(Y)),\]
where $\chi$ means the Euler characteristic with respect to the mod~$2$
grading. In other words:
\be{description}
\item[$\bullet$]If $\del_0,\del'_0$ are both zero then $\zeta=0$.
\item[$\bullet$]If $\del_0\neq0$ then $\zeta=\max\{n\ge1\st\del_{n-1}\neq0\}$.
\item[$\bullet$]If $\del'_0\neq0$ then
$\zeta=-\max\{n\ge1\st\del'_{n-1}\neq0\}$.
\end{description}
\end{defn}

\be{lemma}
$\zeta(-Y,-m)=-\zeta(Y,m)$.
\end{lemma}

\proof The vector space
$\hf^*(-Y,-m)$ is the dual of $\hf_*(-Y,-m)=\hf^{-1-*}(Y,m)$.
Moreover, the $u$--map on the former vector space (in the degrees in which
it is defined) is the dual of the $u$--map on the latter. Similarly,
$\del_0:\hf^{-2m-2}(-Y,-m)\to\f$ is the dual of
$\del'_0:\f\to\hf^{2m+1}(Y,m)$. Together with Proposition~\ref{prop:delnun},
this proves the lemma.\square

\be{prop}\label{basic-ineq}
Let $W$ be a compact, connected $\spc$ $4$--manifold
whose boundary is a disjoint
union of rational homology spheres $Y_1,\dots,Y_r$, where $r\ge1$,
and let $m_j\in\mm(Y_j)$.
Suppose $b^+(W)=0$ and $\zeta(Y_j,m_j)\le0$ for $j>1$. Then
\[\zeta(Y_1,m_1)\ge\sum_{j=1}^rm_j+\frac18(c_1(\cl_W)^2-\si(W))\]
provided the right hand side of the inequality is positive.
\end{prop}
\proof Let $\ga_1,\dots,\ga_\ell$ be a collection of disjoint loops in $W$
representing a basis for $H_1(W;\z)/\text{torsion}$, and let $W'$ denote
the manifold obtained by performing surgery on each $\ga_j$.
Then $b_1(W')=0$. Let $V_i$ be a small open tubular neighbourhood of $\ga_i$,
and set $W^-=W\setminus\cup_iV_i$. Then with integral coefficients
there are isomorphisms
\[H^2(W)\oset\approx\to H^2(W^-)\oset\approx\leftarrow H^2(W'),\]
and similarly for $H_2$. Hence, $W$ and $W'$ have the same
intersection forms, the induced $\spc$ structure on
$W^-$ has a unique extension to $W'$, and
the canonical classes for $W$ and $W'$ have the same square.
Now let $X:=\hw$ be the result of adding half-infinite tubes
to $W$.

Set
\[\vec\theta:=(\theta_{m_1},\dots,\theta_{m_r}).\]
Then
\[\dim M(X;\vec\theta)=2k+1,\]
where
\[k+1=\ind_{\co}(D_{A\rr})=\sum_jm_j+\frac18(c_1(\cl_W)^2-\si(W)).\]
Here $A\rr$ can be taken to be the reference $\spc$ connection over $X$ used
in defining the moduli space $M(X;\vec\theta)$, see \cite[Section~3.4]{Fr13}.

Suppose $k\ge0$. By Proposition~\ref{prop:pert-red} we may assume that
$M(X;\vec\theta)$ contains a unique reducible point,
and that this point is regular.
Let $M^-(X;\vec\theta)$ be the result of removing from
$M(X;\vec\theta)$ a small neighbourhood of the reducible point, so that the
boundary of $M^-(X;\vec\theta)$ is diffeomorphic to $\cp k$.

Let $\bl\to\cb^*(W)$ be the natural
complex line bundle. Choose a generic smooth section $s$ of $\oplus^k\bl$.
For any moduli space $M$ over $X$ and subset $\Si$ of the
irreducible part of $M$ set
\[c_1(\bl)^k\cap \Si:=\{\om\in \Si\st s(\om|_W)=0\}.\]

For $j>1$ the assumption $\zeta(Y_j,m_j)\le0$ means, by the universal
coefficient theorem, that the
cycle $\del\in\cf_{2m_j-2}(Y_j,m_j)$ is in fact a boundary, i.e.\ there is
a $\rho_j\in\cf_{2m_j-1}(Y_j,m_j)$ with
\[\del+\prtl\rho_j=0,\]
where $\prtl$ is the boundary operator of the chain complex $\cf_*$.
For any $\al\in\fl_{Y_1}$ set
\[M_\al:=M(X;\al,\si_2,\dots,\si_r),\]
where $\si_j=\rho_j+\theta_{m_j}$. Here we treat $M(X;\cdot)$ as a multilinear
function, so that $M_\al$ is a formal linear combination of moduli spaces
over $X$. In the case $\al=\theta:=\theta_{m_1}$
one of the terms appearing in that
formal linear combination is $M(X;\vec\theta)$. Let $M^-_\theta$ be obtained
from $M_\theta$ by replacing the term $M(X;\vec\theta)$ with
$M^-(X;\vec\theta)$.

For any irreducible $\al$ set
\[\hm_\al:=c_1(\bl)^k\cap M_\al,\qquad\hm^-_\theta:=c_1(\bl)^k\cap M^-_\theta,\]
where $c_1(\bl)^k\cap\cdot$ is treated as a linear function.
For the formal linear combination $\hm^-_\theta$ of oriented $1$--manifolds the
number of boundary points plus the number of ends (both counted with sign) 
must be zero. The number of boundary points is
\[\#\left(c_1(\bl)^k\cap\prtl M^-(X;\vec\theta)\right)=\pm1.\]
The ends of $\hm^-_\theta$ correspond to factorizations over the ends of $X$.
Because of the choice of $\si_j$, only factorizations over the end
$\rpy_1$ contribute. Therefore, the number of ends of $\hm^-_\theta$ is
$\del a$, where the cocycle $a$ is defined by
\[a=\sum_\al(\#\hm_\al)\cdot\al,\]
the sum being taken over all $\al\in\fl^*_{Y_1}$ of index $2m_1-2$.
Writing $[\ ]$ for cohomology classes we have $[a]=u^k[b]$, where 
\[b=\sum_\beta(\# M_\beta)\cdot\beta.\]
Thus,
\[\del_0u^k[b]=\del\al=\mp1,\]
proving that $\zeta(Y_1,m_1)\ge k+1$ as claimed.\square

\be{cor}\label{cor:smnull}
Let $W$ be a compact, connected $\spc$
$4$--manifold whose boundary is a disjoint
union of rational homology spheres $Y_1,\dots,Y_r$, and let $m_j\in\mm(Y_j)$.
Suppose $b_2(W)=0$ and $\zeta(Y_j,m_j)=0$ for all $j$. Then
$\sum_{j=1}^rm_j=0$.
\end{cor}

\proof If $\sum_jm_j\neq0$ then after perhaps reversing the orientation of
$W$ we may assume $\sum_jm_j>0$, contradicting the proposition.\square

\be{lemma}\label{zeta-s3}
$\zeta(S^3,m)=m$ for all $m\in\z$.
\end{lemma}

\proof Since $\zeta(S^3,-m)=-\zeta(S^3,m)$, we may assume $m\ge0$. The lemma
clearly holds for $m=0$, because $\hf^*(S^3,0)=0$. Now let $m\ge1$. Applying
Proposition~\ref{basic-ineq} with $W=[0,1]\times S^3$, $m_1=m$,
$m_2=0$ we obtain $\zeta(S^3,m)\ge m$.
On the other hand, Proposition~\ref{s3vanish}
gives $\zeta(S^3,m)\le m$. Together this proves the lemma.\square

\be{prop}\label{prop:yyy}
Let $W$ be a $\spc$, compact, connected $4$--manifold with boundary components
$-Y_0,-Y_1,Y_2$, such that each $Y_j$ is a rational homology sphere and
$b_2(W)=0$. Suppose $m_j\in\mm(Y_j)$, $j=0,1,2$ satisfy
$-m_0-m_1+m_2=0$. Then the following hold:
\be{itemize}
\item[(i)]If $\zeta(Y_j,m_j)>0$ for $j=0,1$ then
\[\zeta(Y_2,m_2)\ge\zeta(Y_0,m_0)+\zeta(Y_1,m_1).\]
\item[(ii)]If $\zeta(Y_0,m_0)=0$ then
\[\zeta(Y_1,m_1)=\zeta(Y_2,m_2).\]
\end{itemize}
\end{prop}
\proof The proofs are similar to those of \cite[Lemma~6]{Fr3}
and \cite[Theorem~14]{Fr3}, respectively.\square


\be{cor}\label{zeta-cancel}
If $\zeta(Y_1,m_1)=\zeta(Y_2,m_2)\neq0$ then
$\zeta(Y_1\#(-Y_2),m_1-m_2)=0$.
\end{cor}

\proof Apply part~(i) of the proposition to $W$ or $-W$, where $W$ is the
standard cobordism from $Y_1\cup(-Y_2)$ to $Y_1\#(-Y_2)$.\square

\be{prop}\label{m-exist}
For any $Y$ there exists an $m\in\mm(Y)$ such that $\zeta(Y,m)=0$.
\end{prop}
\proof Pick any $m'\in\mm(Y)$ and set $k:=\zeta(Y,m')$. If $k=0$ then we are
done. Otherwise, set $m=m'-k$.
Then Lemma~\ref{zeta-s3} and Corollary~\ref{zeta-cancel} give
\[\zeta(Y,m)=\zeta(Y\#(-S^3),m'-k)=0.\square\]

\be{prop}\label{zeta-add}
If $\zeta(Y,m)=0$ then $\zeta(Y,m+k)=k$ for all integers $k$.
\end{prop}
\proof Let $D\subset(0,1)\times Y$ be an embedded closed $4$--ball.
Applying Proposition~\ref{prop:yyy}~(ii)
to $W=([0,1]\times Y)\setminus\text{int}(D)$ we get
\[\zeta(Y,m+k)=\zeta(S^3,k)=k.\square\]

\be{defn}For any $\spc$ rational homology $3$--sphere $Y$
we define $h(Y)=h(Y;\f)$ to be the unique $m\in\mm(Y)$ such that $\zeta(Y,m)=0$.
\end{defn}
It follows from Propositions~\ref{m-exist} and \ref{zeta-add}
that $h(Y)$ is well defined, and that
\be{equation}\label{eqn:hmz}
h(Y)=m-\zeta(Y,m)
\end{equation}
for all $m\in\mm(Y)$.

Theorem~\ref{hp-char} is now an immediate consequence
of Propositions \ref{prop:varchamb}, \ref{m-exist}, and \ref{zeta-add},
and Theorem~\ref{thm:h-D} follows from \Ref{eqn:hmz}.

{\em Proof of Theorem~\ref{thm:h-add}:} Set $m_j:=h(Y_j)$
for $j=1,2$ and $m_3:=h(Y_1\#Y_2)$. Applying
Corollary~\ref{cor:smnull} to the standard cobordism $W$ from $Y_1\cup Y_2$ to
$Y_1\#Y_2$ we obtain $-m_1-m_2+m_3=0$, which is the assertion of the theorem.
\square

{\em Proof of Theorem~\ref{main-ineq}:} Set $m_j:=h(Y_j)$ for $j>1$ and
choose $m_1\in\mm(Y_1)$ with
\[\sum_{j=1}^rm_j+\frac18(c_1(\cl_W)^2-\si(W))>0.\]
The theorem then follows from Proposition~\ref{basic-ineq}.\square

We now address the dependence of $h$ on the coefficient field $\f$.
Let $\f$ have characteristic $p$ and let $\bk$ be the prime field of $\f$,
i.e.\ $\bk=\q$ if $p=0$ and $\bk=\z/p\z$ if $p>0$. Then
\[\hf^*(Y;\f)=\hf^*(Y;\bk)\otimes_{\bk}\f,\]
from which one easily deduces that
\[h(Y;\f)=h(Y;\bk).\]
Thus we may write $h_p(Y)$ instead of $h(Y;\f)$.

Note that if $Y$ admits a metric $g$ of positive scalar curvature
then $(g,0)$ is an ND-pair for $Y$ for which the corresponding Floer
cochain complex vanishes. Hence $h_p(Y)=I(g,0)$ for all $p$. (The index $I$
was defined in \Ref{eqn:igv}).

If $h_p(Y)\neq h_0(Y)$ for some $p$ then it is not hard to see that
$\hf^*(Y,h_0(Y);\z)$ contains an element of order $p$ (in degree
$2h_0(Y)\pm1$). 
The following proposition shows that this group is in a sense
the universal home for torsion:

\be{prop}\label{prop:Jtorsion}
For any $m,q$ let $T^q_m$ be the torsion subgroup of $\hf^q(Y,m;\z)$.
For $m_1\le m_2$ the canonical homomorphism
\[\hf^q(Y,m_1;\z)\to\hf^q(Y,m_2;\z)\]
restricts to a homomorphism
\[T^q_{m_1}\to T^q_{m_2}\]
which is injective if $m_2\le h_0(Y)$ and surjective if $h_0(Y)\le m_1$.
\end{prop}
\proof The element $\del'_j\in\hf^*(Y,m_1;\z)$ has infinite order for
$j=0,\dots,h_0(Y)-m_1-1$ if $m_2\le h_0(Y)$ and finite order for all $j$
if $h_0(Y)\le m_1$. The proposition now follows from 
Proposition~\ref{prop:varchamb}.\square


\section{Proof of Theorem~\ref{thm:Lhpwelldef}}
\label{sec:b1b2inv}


\be{lemma}
Let $A,B$ be $r\times r$ complex matrices and $m$ a natural number such that
$\tr(A^n)=\tr(B^n)$ for all natural numbers $n$ satisfying
$m\le n<2r+m$. Then $A$ and $B$ have the same 
characteristic polynomial. In particular, $\tr(A)=\tr(B)$.
\end{lemma}

\proof Let $E$ denote the set of all the eigenvalues of $A$ and $B$.
For any square matrix $C$ and complex number $\lla$
let $m_C(\lla)$ denote the multiplicity of $\lla$ as an eigenvalue of $C$.
It suffices to show that $m_A(\lla)=m_B(\lla)$
for all $\lla\neq0$. So suppose $m_A(\lla)>0$ and consider the polynomial
\[P(z)=z^m\prod_{a\in E\setminus\{\lla\}}(z-a).\]
Putting $A$ in triangular form we see that the multiplicity of $P(\lla)$ as an 
eigenvalue of $P(A)$ is $m_A(\lla)$,
and there are no other non-zero eigenvalues.
Similarly for $B$. Hence
\[m_A(\lla)\cdot P(\lla)=\tr(P(A))=\tr(P(B))=m_B(\lla)\cdot P(\lla).\]
But $P(\lla)\neq0$, so $m_A(\lla)=m_B(\lla)$.\square

\be{cor}\label{cor:v0v1}
Let $U=U_0\oplus U_1$ and $V=V_0\oplus V_1$ be finite-dimensional
mod~$2$ graded complex vector spaces. Let $f=f_0\oplus f_1\in\text{End}(U)$
and $g=g_0\oplus g_1\in\text{End}(V)$ preserve degrees.
If $m$ is a natural number
such that the Lefschetz numbers $\Ltr(f^n)=\Ltr(g^n)$ for all $n\ge m$, then
$\Ltr(f)=\Ltr(g)$.
\end{cor}

\proof Apply the lemma to $f_0\oplus g_1$ and $g_0\oplus f_1$.\square

{\em Proof of Theorem~\ref{thm:Lhpwelldef}:}
Because $Y_j$ is non-separating, there is a loop in $X$ whose intersection
number with $Y_j$ is $\pm1$. Hence each of $Y_0$ and $Y_1$ represents a
generator of $H_3(X;\z)$, and by assumption these two generators are equal.
Let $f_j:X\to S^1$ be a smooth
map such that $1$ is a regular value of $f_j$ and $(f_j)\inv(1)=Y_j$ as
oriented manifolds. Then
$f_0$ and $f_1$ represent the same generator of $[X,S^1]\approx H^1(X;\z)$.
Set
\[X\jinf=\{(x,t)\in X\times\R\st f_j(x)=e^{2\pi it}\}.\]
Then the projection
\[\pi\jinf:X\jinf\to X\]
is a covering whose group of deck transformations is canonically
isomorphic to $\z$. Let 
\[q_j:X\jinf\to\R\]
be the projection onto the second factor. We may take
\[W_j:=(q_j)\inv[0,1]\]
as the definition of $W_j$. Then
\[W\jn:=(q_j)\inv\left[-\left\lfloor n/2\right\rfloor,
n-\left\lfloor n/2\right\rfloor\right]\]
is the result of gluing together $n$ copies of $W_j$ in a chain.

Since $f_0$ and $f_1$ are homotopic,
the coverings $\pi\zinf$ and $\pi_{1,\infty}$
are isomorphic. In particular, $\pi\zinf$ restricts to a trivial covering
of $Y_1$. Let $B$ be a component of $(\pi\zinf)\inv(Y_1)$. Since $B$ is
compact, there is a natural number $m$ such that $B$ is contained in the
interior of $W_{0,m}$. For every natural number $n$ set
$X_n:=X\zinf/n\z$ and let
\[\pi_n:X_n\to X\]
be the natural $n$--fold covering.
For $n\ge m$ there is a component $C_j$ of $\pi_n\inv(Y_j)$, $j=0,1$, with 
$C_0\cap C_1=\emptyset$. Let $C_j$ be oriented such that the diffeomorphism
$C_j\to Y_j$ obtained by restriction of $\pi_n$ is orientation preserving.
Then $X_n$ is the union of two compact, connected,
codimension~$0$ submanifolds $Z_0,Z_1$ with 
\[Z_0\cap Z_1=C_0\cup C_1,\qquad\prtl Z_j=(-C_j)\cup C_{1-j}.\]
Moreover, $b_1(Z_j)=0=b^+(Z_j)$. If each $Y_j$ is an integral homology sphere
then it follows from Theorem~\ref{main-ineq} and Elkies' theorem
that $h_p(Y_0)=h_p(Y_1)$. The same conclusion holds if $b_2(X)=0$,
because then $b_2(Z_j)=0$, too.

Returning to the general case, we can identify
\[W\jn=Z_{1-j}\cup_{C_{1-j}}Z_j.\]
Setting
\[A_j:=\rpsi(W_j):\rhf^*(Y_j)\to\rhf^*(Y_j)\]
we have
\[(A_j)^n=\rpsi(\ww jn)=\rpsi(Z_{1-j})\circ\rpsi(Z_j).\]
Since $\rpsi(Z_j)$ preserves the mod~$2$ grading, we conclude that
\[\Ltr((A_0)^n)=\Ltr((A_1)^n)\quad\text{for all $n\ge m$.}\]
Corollary~\ref{cor:v0v1} now says that $\Ltr(A_0)=\Ltr(A_1)$,
and the theorem is proved.
\square

\section{A finite-dimensional analogue}
\label{sec:fd-analogue}

Whereas Marcolli--Wang's construction of equivariant Seiberg--Witten
Floer homology \cite{Marcolli-Wang1} was modelled on
Austin--Braam's equivariant Morse complex \cite{AuB1},
it may seem less obvious why the groups $\ohf^*(Y)$
deserve being called {\em equivariant} Floer cohomology. 
This section will provide justification for this name by means of a
finite-dimensional analogue. We will use integer coefficients, for 
simplicity.

Conceptually, a natural setting for our construction of Floer cohomology
is the Hilbert manifold $\cb_y(Y)$ of
$L^2_1$ configurations $(B,\Psi)$ on $Y$ modulo $L^2_2$ gauge transformations
$u:Y\to\U1$ with $u(y)=1$, where $y\in Y$ is a base-point.
There is a smooth action of $\U1$ on $\cb_y(Y)$ whose fixed-points
are the elements $[B,\Psi]$ with spinor $\Psi=0$.
Away from the fixed-point set $F$
the action is free, so that the quotient
\[\cb^*(Y)=(\cb_y(Y)\setmin F)/\U1\]
is also a Hilbert manifold.  The Chern--Simons--Dirac functional on $Y$
descends to a smooth $\U1$--invariant function
\[\ucd:\cb_y(Y)\to\R,\]
which in turn induces a smooth function
\[\uucd:\cb^*(Y)\to\R.\]
An important point here is that $\ucd$ has exactly one critical point
which is also a fixed-point.
The group $\hf^*(Y,m)$ was defined as the cohomology of a ``Morse complex''
of $\uucd$, whereas $\ohf^q(Y)$ was obtained from $\hf^*(Y,m)$ by 
a limiting construction as $m=\frac12\ind(\theta)\to\infty$.

We will now interpret this construction in terms of finite-dimensional
Morse theory (see for instance \cite{Banyaga-Hurtubise}).
As our analogue of $\cb_y(Y)$ we consider a compact
manifold-with-boundary $V$ of positive dimension
on which $\U1$ acts smoothly, and
such that the action is free away 
from one fixed-point $p$ which lies in the interior of $V$.
In the proof of Proposition~\ref{prop:findim} below we will show that
the $\U1$--action on $V$ around $p$ is equivalent to the standard action on 
$\co^n$ around the origin for some $n$.
In particular, $\dim V=2n$. (It is easy to
deduce from this that the boundary $\prtl V$ must be non-empty.) 
We assume
$\prtl V$ is the disjoint union of two relatively open subsets $V_0$ and $V_1$,
one of which might be empty. Let
\[\pi:V\to B:=V/\U1\]
be the projection onto the quotient space
and set $\ip:=\pi(p)$ and $B_i:=\pi(V_i)$ for $i=0,1$. Then the manifold
\[B^*:=B\setmin\{\ip\}.\]
is our analogue of $\cb^*(Y)$. Now let 
\[f:V\to[0,1]\]
be a smooth $\U1$--invariant function such that the following hold:
\be{itemize}
\item $f\inv(i)=V_i$ for $i=0,1$.
\item $p$ is a non-degenerate critical point of $f$ of maximal index,
i.e.\ $\ind(p)=2n$. (In particular, $p$ is a local maximum of $f$.)
\item The induced function
\[g:B^*\to\R\]
has only non-degenerate critical points, all of which lie in the interior
of $B^*$.
\end{itemize}
Note that taking $\ind(p)$ maximal mirrors the construction of $\ohf^*(Y)$.

Choose a $\U1$--invariant Riemannian metric on $V$
and let $B^*$ have the Riemannian metric for which 
\[V^*:=V\setmin\{p\}\to B^*\]
is a Riemannian submersion. Then flow-lines of $-\nabla f$ in $V^*$
map to flow-lines of $-\nabla g$ in $B^*$.
Let $C$ be the image in $B$ of the set of 
critical points of $f$. Then
\[C^*:=C\setmin\{\ip\}\]
is the set of critical points of $g$.
For any $a\in C$ let $W^s(a)$ and $W^u(a)$ denote the stable and unstable
manifolds of $a$, resp. (The usual definitions carry over in an obvious way,
even though $B$ is not a manifold due to the singular
point $\ip$.) In particular, $W^u(\ip)$ agrees with the image in $B$
of the unstable manifold of $p$ in $V$, and $W^s(\ip)=\{\ip\}$.
We assume that the Morse--Smale
condition holds in the sense that $W^s(a)$ and $W^u(b)$ intersect
transversely for any $a\in C^*$, $b\in C$ (which holds trivially if $b=\ip$),
and set
\[M(a,b):=W^s(a)\cap W^u(b).\]
Then
\[\dim M(a,\ip)=2n-1-\ind(a).\]
For comparison: In gauge theory one has, for any $\al\in\fl^*_Y$,
\[\dim M(\al,\theta)=2m-1-\ind(\al).\]
This means that $2n=\ind(p)$ corresponds to $2m=\ind(\theta)$.
The analogue of $\hf^*(Y,m)$ is now the cohomology of the usual Morse
complex of $g$, which is isomorphic to the singular cohomology of the pair
$(B^*,B_0)$.

Before stating the main result of this section recall that
if a Lie group $K$ acts on a space $X$ then the {\em equivariant cohomology}
$H^*_K(X)$ is by definition the cohomology of the homotopy quotient
\[X_K:=X\times_K EK,\]
where $EK\to BK$ is a universal principal $K$--bundle.
If $Y$ is a $K$--invariant subspace of $X$ such that $Y_K\subset X_K$ has
the subspace topology (which it always has when $K$ is compact and $Y$ is 
closed in $X$)
then we define
\[H^*_K(X,Y):=H^*(X_K,Y_K).\]
If the projection $X\to X/K$ is a principal bundle then
the canonical map $X_K\to X/K$ is a weak
homotopy equivalence and therefore induces an isomorphism on (co)homology.
For instance, this applies to the map
\[(V^*)_{\U1}\to V^*/\U1=B^*.\]
Therefore,
the inclusions of $(V^*,V_0)$ and $(\{p\},\emptyset)$ in $(V,V_0)$
induce homomorphisms
\be{equation}\label{eqn:hqb}
H^q(B^*,B_0)\leftarrow H^q_{\U1}(V,V_0)\rightarrow H^q(\cp\infty).
\end{equation}

\be{prop}\label{prop:findim}
The first map in \Ref{eqn:hqb} is an isomorphism for $q\le2n-1$ whereas the
second map is an isomorphism for $q\ge2n$.
\end{prop}
Thus, the equivariant cohomology of $(V,V_0)$ (as an Abelian group)
is completely determined by the cohomology of $(B^*,B_0)$. 
An analogous statement holds for the Floer groups provided
$m$ is so large that $\uhf^q(Y)=0$ for $q\ge2m$. Namely, in the
sequence
\be{equation}\label{eqn:gaugeth-maps}
\hf^q(Y,m)\to\ohf^q(Y)\oset D\to P^q(Y)
\end{equation}
the first map 
is an isomorphism for $q\le2m-1$ whereas $D$ is an isomorphism
for $q\ge2m$. (See Subsection~\ref{subsec:equiv-floergr}.)
Thus, in this analogy with finite-dimensional Morse theory, $\ohf^*(Y)$
plays the role of the $\U1$--equivariant cohomology of $\cb_y(Y)$
(since the latter has no boundary).


{\em Proof of Proposition~\ref{prop:findim}:} Set $\bt:=\U1$.
Note that the action of $\bt$ on $V$ induces a representation of $\bt$ on the
tangent space $T_pV$. Let $D'\subset T_pV$ be a small closed ball centred at the
origin. Then $\exp:D'\to V$ is a $\bt$--equivariant embedding.
Because $\bt$ acts freely on 
$V^*$ it follows from the classification of real
$\bt$--representations that for some integer $n$ there is an isomorphism
$T_pV\oset\approx\to\co^n$ of euclidean vector spaces which is equivariant with
respect to the standard action of $\bt$ on $\co^n$.

Let $S$ be the boundary of the geodesic ball $D:=\exp(D')$.
Note that there is a deformation retraction of the homotopy quotient $D_\bt$
onto a subspace which we can identify with $B\bt=\cp\infty$. Moreover, the
inclusion $S_\bt\to D_\bt$ is homotopic to a classifying map
$S_\bt\to\cp\infty$ of the $\bt$--bundle $S\times E\bt\to S_\bt$.
From the Mayer--Vietoris sequence of the pairs $((V^*)_\bt,(V_0)_\bt)$
and $(D_\bt,\emptyset)$ in $(V_\bt,(V_0)_\bt)$ we
then obtain for every integer $q$ a short exact sequence
\[0\longrightarrow H^q_\bt(V,V_0)
\longrightarrow H^q(B^*,B_0)\oplus H^q(\cp\infty)\oset{a+b}
\longrightarrow H^q(\cp{n-1})\longrightarrow0,\]
where $b$ is an isomorphism for $q\le2n-1$. Since
\[H^q(B^*,B_0)=H^q(\cp{n-1})=0\quad\text{for $q\ge2n$}\]
for dimensional reasons, the proposition follows.
\square

\noindent\textsc{Institut f\"ur Mathematik,\\
Universit\"at Z\"urich,\\
Winterthurerstrasse 190,\\
CH-8057 Z\"urich,\\
Switzerland}
\\ \\
E-mail:\ froyshov@math.uzh.ch

\end{document}